\documentclass[10pt]{article}
\usepackage{amsmath, amsfonts, amsthm}
\usepackage{amssymb,mathrsfs, verbatim}
\usepackage{graphicx,amstext}
\usepackage{amsmath}
\usepackage{mathtools}
\usepackage{amssymb}
\usepackage{latexsym}
\usepackage{ mathrsfs }
\usepackage{enumerate}
\usepackage[USenglish]{babel}

\usepackage[latin1]{inputenc}
\newcommand{\clg}[1]{{\mathcal{#1}}}

\newcommand{\R}{\mathbb R}

\newcommand{\vp}{\varphi}

\newcommand{\sgn}{\text{sgn}}

\newcommand \loc    {\text{loc}}
\newcommand{\ess}{\,{\rm ess}}

\newcommand{\dive}{{\rm div}}

\numberwithin{equation}{section}

\newtheorem{theorem}{Theorem}[section]
\newtheorem{proposition}{Proposition}[section]
\newtheorem{remark}{Remark}[section]
\newtheorem{lemma}{Lemma}[section]
\newtheorem{corollary}{Corollary}[section]
\newtheorem{definition}{Definition}[section]

\newcommand{\intav}[1]{\mathchoice {\mathop{\vrule width 6pt height 3 pt depth  -2.5pt
\kern -8pt \intop}\nolimits_{\kern -6pt#1}} {\mathop{\vrule width
5pt height 3  pt depth -2.6pt \kern -6pt \intop}\nolimits_{#1}}
{\mathop{\vrule width 5pt height 3 pt depth -2.6pt \kern -6pt
\intop}\nolimits_{#1}} {\mathop{\vrule width 5pt height 3 pt depth
-2.6pt \kern -6pt \intop}\nolimits_{#1}}}

\begin{document}
\pagestyle{myheadings} \markboth{\textsc{G. Huaroto$\&amp;$ W.
Neves}}{\textsc{Fractional Chemotaxis Models}}
\title{Solvability of the Fractional Hyperbolic Keller-Segel System}

\author{Gerardo Huaroto and Wladimir Neves}

\date{}
%%%%%%%%%%%%%%%%%%%%
\maketitle

%%%%%%%%%%%%%%%%%%%%%%%%%%%%%%%%%
\begin{abstract}
We study a new nonlocal approach to the mathematical modelling of the 
Chemotaxis problem, which describes the random motion of a certain population 
due a substance concentration.
Considering the initial-boundary value problem for the fractional
hyperbolic Keller-Segel model, we prove the solvability of the problem. 
The solvability result relies mostly on 
fractional calculus and kinetic formulation of scalar conservation laws.
%the kinetic method.
\end{abstract}
%%%%%%%%%%%%%%%%%%%%%%%%%%%%%%%%%

%%%%%%%%%%%%%%%%%%%%
\section{Introduction} 
\label{INTROD}
%%%%%%%%%%%%%%%%%%%%
We introduce and study in this paper the  Fractional Hyperbolic Keller-Segel (FHKS for short) model for chemotaxis
described by the following system
\begin{equation}
\label{eq.main}
\left \{
\begin{aligned}
    &\partial_t u + \dive \big(g(u) \,  \nabla \mathcal{K}_s c \big)= 0, \quad  \text{in $(0,\infty) \times \Omega$},
\\[5pt]
    &  (-\Delta_N)^{1-s} \ c + c= u,  \quad \text{in $\Omega$},
\\[5pt]
    & u|_{\{t=0\}}= u_{0},  \hspace{21pt} \text{in $\Omega$},
\\[5pt]
 & \nabla\mathcal{K}_s c\cdot \nu= 0, \hspace{22pt}  \text{on $\Gamma$},
\end{aligned}
\right.
\end{equation}
where $u(t,x)$ is the density of cells and $c(t,x)$ is the
chemoattractant concentration, which is responsible for the cell
aggregation. The problem is 
posed in a bounded open subset $\Omega \subset \mathbb R^n$, ($n=1,\,2$, or 3),
with $C^2-$boundary denoted by $\Gamma$, and as usual we denote by 
$\nu(r)$ the outward normal to $\Omega$ at $r \in \Gamma$. 
The given measurable bounded function $u_{0}$ is the initial  
condition of the cells, and we assume 
\begin{equation}
\label{ID}
\text{$0 \leq u_0(x) \leq 1$, \; for a.e. $x \in \Omega$} .
\end{equation}
Moreover, since the normal fractional flux in the equation on $u$ vanishes on $\Gamma$,
that is, we assume that the boundary is characteristic, it is not necessary to prescribe 
boundary conditions for $u$. This assumption besides the natural one, it prevents 
some specific difficulties related to the trace problem, see \cite{WNEPJS} for instance. 

\medskip
Here for $0<s< 1$, $(-\Delta_N)^s$ denotes the Neumann spectral fractional Laplacian (NSFL for short)
operator, which characterizes long-range diffusion effects.
We also consider the non-local operator $\mathcal{K}_s$ given by definition in \eqref{DEFFRALAPINV}, 
which can not be interpreted as the inverse of $(-\Delta_N)^s$. Indeed, we are dealing with the spectral fractional Laplacian
$(-\Delta_N)^s$, which is
defined from the eigenvalues of the Laplacian with Neumann boundary conditions, where
the first eigenvalue is zero and the corespondent eigenfunction is a constant. Hence
the NSFL is not an injective operator in its domain, actually it becomes injective when restrict to the set
\begin{equation}
\label{MEDIA0}
   \Big\{ f \in D\big((-\Delta_N)^{s}\big); \, \int_{\Omega} f(x)\,dx= 0 \Big\}. 
\end{equation}
Therefore, the operator $((-\Delta_N)^s)^{-1} =:(-\Delta_N)^{-s}$ does not exist in general, and thus we are not allowed to write, 
$(-\Delta_N)^{1-s}c= (-\Delta_N) (-\Delta_N)^{-s}c$, unless $c(t,x)$ satisfies \eqref{MEDIA0}. 
The reader is addressed to Section \ref{FNL} for a comprehensive description of the NSFL operator.
%and the inverse of the restricted NSFL operator on the above set. 

\medskip
%On the other hand, since $1-s> 0$ the operator $(-\Delta_N)^{1-s}$ is well defined in its domain,
%and 
Now, due to the first eigenfunction of  $(-\Delta_N)^{1-s}$ be $\vp_0= 1/ \sqrt{|\Omega|}$, 
it is not difficult to show that 
\begin{equation}
\label{MEDA}
\int_\Omega (-\Delta_N)^{1-s} \; c(t,x) \, dx= 0.
\end{equation}
Indeed, we have for almost all $t> 0$
$$
   \int_\Omega (-\Delta_N)^{1-s} \; c(t,x) \, dx= \sqrt{ |\Omega|} < (-\Delta_N)^{1-s} c(t,\cdot), \vp_0>= 0, 
$$
where we have used Proposition \ref{THMCARAC} item (2). 
Consequently, from the 
second equation in \eqref{eq.main} it follows that, for a.a. $t> 0$
\begin{equation}
   \int_\Omega c(t,x) \, dx\,=\, \int_\Omega u(t,x) \, dx, \label{eq_contabilidade}
\end{equation}
which implies that $c$ satisfies condition \eqref{MEDIA0} if, and only if, $u$ satisfies it.
%can not satisfy condition \eqref{MEDIA0} unless $u$ satisfies it. 
 
Last but not least, we consider (conveniently) the function $g(u)= u \, (1-u)$, 
which prevents blow up.
Indeed, from a maximum principle established result, 
and the assumption \eqref{ID}
we are going to show that, 
$0 \leq u(t,x) \leq 1$ almost everywhere. Consequently, 
we could not impose that the chemoattractant concentration
$c(t,x)$ satisfies condition \eqref{MEDIA0}, otherwise $u= 0$ a.e.. 
This implies an inherent fractional characterization of the problem studied here. 

\medskip
The theory of chemotaxis modeling  
goes back to E. F. Keller and L. A. Segel \cite{KS0, KS1,KS2}, where a
detailed description of the movement of cells oriented by chemical cues
can be found.  
In fact, a nonlocal version of the Keler-Segel model has been proposed
by Caffarelli, Vazquez in \cite{caff-vazquez-arma} 
(Section 8. Comments and Extensions).
Clearly, the study of a nonlocal version of the Keler-Segel model 
becomes interesting as a proposed open problem in that seminal cited paper 
\cite{caff-vazquez-arma}. But not just because of that, 
it seems natural to assume that the random motion of a certain population 
due a substance concentration can be described by L\'evy process 
instead of the Brownian motion, which is in fact a special
case of the continuous L\'evy process.
The fractional model proposed here in \eqref{eq.main}
is a fractional generalization of the model considered in 
Perthame-Dalibard \cite{perth-dali-tams}. Indeed, in that paper 
they studied the following system 
\begin{equation}
\label{DAPER}
\begin{cases}
    \partial_t u+\dive \big(g(u) \,  \nabla S \big)= 0, \hspace{9pt} \text{in $(0,\infty) \times \Omega$},
\\[5pt]
      (-\Delta) S + S= u,  \hspace{5pt} \text{in $\Omega$},
\\[5pt]
     u|_{\{t=0\}}= u_{0},  \hspace{19pt} \text{in $\Omega$},
\\[5pt]
 \nabla S \cdot \nu= 0, \hspace{25pt} \text{on $\Gamma$},
\end{cases}
\end{equation}
which follows from the system \eqref{eq.main}, at least formally
passing to the limit as $s \to 0^+$. In particular, it is allowed to have jumps  
in the chemoattractant concentration
$c(t,x)$ in the proposed model 
\eqref{eq.main}, which is not the case for the original model \eqref{DAPER}.
This is very important for practical applications. 

\medskip
The mathematical analysis developed in this paper combine two different 
aspects (fractional calculus and kinetic formulation of scalar conservation laws).
The first one is quite new, where well-posedness and several useful estimates
are obtained for a fractional parabolic-elliptic system, see Section \ref{ParabolicApp}. 
Then, the rigidity result of Perthame-Dalibard 
\cite{perth-dali-tams} is well adapted in Section \ref{LP} to tackle 
the limiting problem to obtain the proof of the Theorem \ref{ETGS} (Main Theorem). 
Indeed, the non-local description of the chemotaxis model studied here 
brings new difficulties to the entropy structure of the first equation in \eqref{eq.main}, 
see \eqref{FORMAL}.  
In particular, because of the non-local space dependence of the entropy flux,
which leads to specific difficulties.  
More precisely, due to the lack of strong compactness or contraction properties, 
we apply the kinetic formulation associated to the
propagation of oscillations of the transport equations (Theorem \ref{FF}). To this end, it is essential to use
the renormalization procedure in the kinetic equation,
and we are allowed to do that for any $s \in (0,1/2]$. 
%The solvability 
%of \eqref{eq.main} for $s \in (1/2, 1)$ remains open,
%which could not be resolved by the technics used here. 

\medskip
We would like to remark that the uniqueness property 
is not addressed in this paper, 
neither the asymptotic behavior. 
First, let us stress that no uniqueness result is 
expected 
%has been proven 
even for the standard system \eqref{DAPER}. On the other hand, 
long-time behavior was establish to \eqref{DAPER} in
\cite{perth-dali-tams}, and similarly the main ingredient require here is 
$$
  \int_\Omega g(u(t)) \nabla \mathcal{K}_s c(t) \cdot \nabla c(t,x) \, dx \geq 0
$$
for each $t> 0$, 
which does not seen that such 
inequality holds. Indeed, we may write 
$$
\begin{aligned}
  \int_\Omega g(u(t)) \nabla \mathcal{K}_s c(t) \cdot \nabla c(t) \, dx 
 &= \int_\Omega g(c(t)) \nabla \mathcal{K}_s c(t) \cdot \nabla c(t) \, dx
 \\[5pt]
 & +  \int_\Omega (1-2 \xi) (u(t) - c(t)) \nabla \mathcal{K}_s c(t) \cdot \nabla c(t) \, dx,
\end{aligned}
$$
where $\xi \in (\min\{u,c\}, \max\{u,c\})$. Formally, the first integral in the right hand side is positive
(see Caffarelli, Soria, V\'azquez
\cite{caff-vazquez-soria-jems} p.1706), although we do not have sign control in the second one. 
We leave this question for future studies.

\medskip
Finally, we discuss in Section \ref{Comments} some related problems that can be handled similarly. 
In particular, one of the systems considered does not satisfy exactly the condition \eqref{eq_contabilidade},
compare it with \eqref{NEWCOND}. 
%Moreover, we highlight that
%the theory of fractional evolutionary partial 
%differential equations posed in bounded domains
%is an interesting area. Indeed, 
%One of the major reasons is that, 
%there exist different non-equivalent definitions of the 
%there is not a unified way to define 
%fractional Laplacian operator in bounded domains,
%and each specific problem may have a correct choice.
%in particular with respect to boundary conditions. 
%For instance, we address the reader to the authors papers \cite{GHWN1}, and 
%\cite{GHWN2}. 

%%%%%%%%%%%%%%%%%%%%%%
\subsection{Statement of the FHKS system}
\label{SFHKSM}
%%%%%%%%%%%%%%%%%%%%%%

The aim of this section is to formulate the mathematical problem for the FHKS system.
We begin observing that, 
the first equation in \eqref{eq.main} is a hyperbolic scalar conservation law, thus the
density of cells function $u$ may admit shocks. Therefore, in order to select the
more correct physical solution, we need an admissible criteria, which is given by the entropy condition. 
%as given at the following

\begin{definition}
\label{DEFENT}
A pair $\mathrm{F}(u)= (\eta(u), q(u))$ is called an entropy pair for the first equation in \eqref{eq.main}, 
if there exists $\eta: \mathbb{R} \rightarrow \mathbb{R}$ a Lipschitz continuous and also convex 
function and the function $q: \mathbb{R} \rightarrow \mathbb{R}$, which satisfies
$$
q^{\prime}(u)=\eta^{\prime}(u) g^{\prime}(u),
$$
for almost all $u \in \mathbb{R}$. Also, we call $\eta(u)$ an entropy and $q(u) \nabla \mathcal{K}_s c$ the associated entropy flux
for the first equation in \eqref{eq.main}. 
\end{definition}
Analogously to the most important example of the Kru\v zkov's entropies, we consider
$$
   \mathbf{F}(u, v)=(|u-v|, \operatorname{sgn}(u-v)(g(u)-g(v)))
$$
for each $v \in \R$.  Another two examples of parameterized family of entropy pairs for \eqref{eq.main},
which will be conveniently used for the kinetic formulation,  
are given by
$$
\mathbf{F}^{\pm}(u, v)=\left(|u-v|^{\pm}, \operatorname{sgn}^{\pm}(u-v)(g(u)-g(v))\right)
$$
for each $v \in \mathbb{R}$, where $|v|^{\pm}:=\max \{\pm v, 0\}$, and
$$
\operatorname{sgn}^{+}(v):=\left\{\begin{array}{ll}
1, & \text { if } v>0 \\
0, & \text { if } v \leqslant 0,
\end{array} \quad \operatorname{sgn}^{-}(v):=\left\{\begin{array}{ll}
\!\! -1, & \text { if } v<0 \\
\; 0, & \text { if } v \geqslant 0.
\end{array}\right.\right.
$$

\medskip
In order to formulate the mathematical problem for the system \eqref{eq.main}, we formally assume enough regularity to 
the pair $(u(t,x), c(t,x))$. Then for any entropy $\eta \in C^2$, we  
multiply by $\eta'(u)$ the fist equation in \eqref{eq.main} to obtain
\begin{equation}
\label{FORMAL}
    \partial_t \eta(u) + {\rm div} \big(q(u) \nabla \mathcal{K}_s c \big)
    + \big(u - c \big) [q - g \eta'](u) \leq 0, 
\end{equation}
where we have used that, $-{\rm div} \nabla \mathcal{K}_s c= u-c$. 
One recalls that, any smooth entropy pair 
$\mathbf{F}(u)= (\eta (u),q(u))$ for \eqref{eq.main} can be recovered by the family 
of Kru\v zkov's entropies. Therefore, the following definition tells us in which sense 
a pair of functions $(u,c)$ is a weak solution of $\textbf{FHKS}$ system: 
\eqref{eq.main}-\eqref{ID}.

\begin{definition}
\label{FHKS} 
Given an initial data $u_0 \in L^\infty(\Omega)$ satisfying \eqref{ID} and any $s \in (0,1)$,  
a pair of functions
\begin{equation*}
(u, c)\in L^{\infty }((0,\infty) \times \Omega) \times L^{\infty}((0,\infty); D((-\Delta_N)^{(1-s)/2}))
\end{equation*}
is called a weak solution to the $\textbf{FHKS}$ system, if for almost all $t> 0$,  the pair $(u,c)$ 
satisfies the condition \eqref{eq_contabilidade}, and the integral inequality
\begin{align}
 \label{CLE}
& \int_0^\infty \!\!\!\! \int_{\Omega} |u(t,x)-v| \;  \partial_t \phi \, dxdt  
+\int_{\Omega} |u_{0}(x)-v| \, \phi (0) \ dx
\notag 
\\[5pt]
& \;\; + \int_0^\infty \!\!\!\! \int_{\Omega} \mathrm{sgn}(u(t,x)-v)
\Big( \big(g(u(t,x))-g(v)\big)\;
\nabla \mathcal{K}_s c(t,x) \cdot \nabla \phi + g(v) \phi \Big) \geq 0 
\end{align}
for any fixed $v\in \mathbb{R}$, and each nonnegative function $\phi \in
C_{c}^{\infty }(\R \times \mathbb{R}^{d})$ and also the following
integral identity
\begin{equation}
\label{EE}
\int_{\Omega } \big(\nabla \mathcal{K}_s c \cdot \nabla \psi + c \ \psi \big) \,dx
= \int_\Omega u \ \psi \ dx,  
\quad \text{for a.a. } t> 0
\end{equation}%
holds for any $\psi \in H^{1}(\Omega )$.
\end{definition}

\medskip
Now, we are able to state plainly the main
result of this paper. Then, we have the following 

\begin{theorem}[Main Theorem]
\label{ETGS} Let $u_0 \in L^\infty(\Omega)$ be an initial data satisfying \eqref{ID} and $0< s \leq 1/2$.
Then, there exists a pair of functions
\begin{equation*}
(u, c)\in L^{\infty }((0,\infty) \times \Omega) \times L^{\infty}((0,\infty); D((-\Delta_N)^{1-s})),
\end{equation*}
which is a weak solution to the $\textbf{FHKS}$ system, and it satisfies 
$$
   0 \leq u(t,x) \leq 1, \quad 0 \leq c(t,x) \leq 1,
$$
for almost all $t> 0$ and $x \in \Omega$. 
\end{theorem}

\begin{remark} 
\label{REM} In fact, the condition \eqref{eq_contabilidade} follows from \eqref{EE} with an integration by parts,
due to the regularity of the function $c(t,\cdot) \in D((-\Delta_N)^{1-s})$ for almost all $t> 0$. Indeed, 
the second equation in \eqref{eq.main} is satisfied almost everywhere, and 
\eqref{eq_contabilidade} follows integrating it on $\Omega$. 
\end{remark}

%%%%%%%%%%%%%%%%%%%%%%%
\subsection{Notation and Functional Spaces}
%%%%%%%%%%%%%%%%%%%%%%

Let $\Omega$ be a bounded open set in $\R^n$.
We denote by $dx$ (or $d\xi$, etc.)
the Lebesgue measure, and by $\clg{H}^\theta$ the $\theta-$dimensional Hausdorff
measure. By $L^p(\Omega)$ we denote the set of real $p-$summable functions 
with respect to the Lebesgue measure,
and the vector counterparts of these spaces are denoted by
$\mathbf{L}^p(\Omega)= \big(L^p(\Omega)\big)^n$.
% the Cartesian product of $L^2(\Omega)$
%$n-$times. 
% (vector ones should be understood 
%componentwise).

\bigskip
%%%%%%%%%%%%%%%%%%%%%%%%%%%%%%%%%%%%%%%%%%%%%%%%%%%%%
$\bullet$ {\bf The space $W^{s,p}(\Omega)$}
%%%%%%%%%%%%%%%%%%%%%%%%%%%%%%%%%%%%%%%%%%%%%%%%%%%%%%%%%%

The fractional Sobolev space is denoted by $W^{s,p}(\Omega)$, where a 
real $s\geqslant 0$ is the
smoothness index, and a real $p\geqslant 1$
is the integrability index. More precisely, for $s \in (0,1)$, 
$p \in [1,+\infty)$, the fractional Sobolev space of
order $s$ with Lebesgue exponent $p$ is defined by
$$
   W^{s,p}(\Omega):= \Big\{ u \in L^p(\Omega): \int_{\Omega} \int_{\Omega}\dfrac{\vert u(x)-u(y)\vert^p}{
\vert x-y\vert ^{n+sp}} \ dx dy< + \infty \Big\},
$$
endowed with norm
$$
   \Vert u\Vert_{W^{s,p}(\Omega)}= \left(  \int_{\Omega}\vert u\vert^p dx + \int_{\Omega} \int_{\Omega}\dfrac{\vert u(x)-u(y)\vert^p}{
    \vert x-y\vert ^{n+sp}}dxdy \right)^\frac{1}{p}.
$$
Moreover, for $s > 1$ we write $s = m + \sigma$, where $m$ is an integer
and $\sigma \in (0, 1$). In this case, the space $W^{s,p}(\Omega)$ consists of those equivalence classes
of functions $u \in W^{m,p}(\Omega)$ whose distributional derivatives $D^{\alpha} u$, with $|\alpha| = m$,
belong to $W^{\sigma,p}(\Omega)$, that is
$$
        W^{s,p}(\Omega)= \Big\{ u \in W^{m,p}(\Omega): 
        {\displaystyle \sum_{\vert \alpha\vert = m}}\Vert D^{\alpha}u \Vert_{W^{\sigma,p}(\Omega)}< \infty \Big\},
$$
which is a Banach space with respect to the norm
$$
\Vert u\Vert_{W^{s,p}(\Omega)}= \Big(  \Vert u\Vert^p_{W^{m,p}(\Omega)} 
+ {\displaystyle \sum_{\vert u\vert = m}}\Vert D^{\alpha}u \Vert^p_{W^{\sigma,p}(\Omega)} \Big)^\frac{1}{p}.
$$
If $s = m$ is an integer, then the space $W^{s,p}(\Omega)$ coincides with the Sobolev space
$W^{m,p}(\Omega)$. Also, it is very interesting the case when $p = 2$, i.e. $W^{s,2}(\Omega)$. 
In this case, the (fractional)
Sobolev space is also a Hilbert space, and we can consider the inner product
$$
    \langle u,v\rangle_{W^{s,2}(\Omega)}= \langle u,v\rangle + \int_{\Omega} \int_{\Omega}
    \frac{(u(x)-u(y))}{\vert x-y\vert ^{\frac{n}{2}+s}} \ \frac{(v(x)-v(y))}{\vert x-y\vert ^{\frac{n}{2}+s}} \ dx dy,
$$
where $\langle\cdot,\cdot\rangle$ is the inner product in $L^2(\Omega)$. 
%On the other hand, we can define the subspace $W_0^{s,p}(\Omega)$ as
%$$
%W^{s,p}_0(\Omega)=\overline{C^{\infty}_c(\Omega)}^{\|\cdot\|_{W^{s,p}(\Omega)}}.
%$$

\bigskip
%%%%%%%%%%%%%%%%%%%%%%
$\bullet$ {\bf The space $H^s(\Omega)$}
%%%%%%%%%%%%%%%%%%%%%%

Following Lions, Magenes \cite{LionsMagenes}, we can define for $s\in (0,1)$, the spaces $H^s(\Omega)$ 
by interpolation between $H^1(\Omega)$ and $L^2(\Omega)$, i.e.
$$
H^s(\Omega)=[H^1(\Omega),L^2(\Omega)]_{1-s}.
$$
According to this definition, this space is a Hilbert space with the natural norm given by the interpolation.
We recall that, when $\Omega$ has Lipschitz boundary regularity, 
then the spaces $W^{s,2}(\Omega)$ and $H^s(\Omega)$ are
equivalent. 

%\begin{theorem}[Trace]
%Let $\Omega$ be a bounded open set in $\R^n$ with Lipschitz boundary. 
%For $1/2 < s \leq 1$, there exists a linear operator 
%$T: H^s(\Omega) \rightarrow L^2(\partial\Omega)$.
%\end{theorem}
%\begin{proof}
%See Lions, Magenes \cite{LionsMagenes} Theorem 9.4, and Theorem 11.5.
%\end{proof}

%%%%%%%%%%%%%%%%%%%%%%%%%%%%%%%
\section{The NSFL Operator}
\label{FNL}
%%%%%%%%%%%%%%%%%%%%%%%%%%%%%%%

Here and subsequently $\Omega \subset \R^n$ is a bounded open set with $C^2-$boundary $\Gamma$. 
Following \cite{GHWN1,GHWN2}, we are interested in fractional powers of a strictly positive self-adjoint 
operators defined in a domain, which is dense in a (separable) Hilbert space. 
More precisely, let us denote by $(-\Delta_{N})$ the operator $(-\Delta)$ subject to Neumann 
boundary conditions. One observes that $(-\Delta_{N})$ is a nonnegative and self-adjoint 
operator defined in 
$$
\begin{aligned}
     D(-\Delta_{N})&= \left\lbrace u\in H^1(\Omega):(-\Delta)u\in L^2(\Omega),\,\mbox{ with } \nabla u\cdot\nu=0 \mbox{ on } \Gamma \right\rbrace
\\[5pt]
    &=\left\lbrace u\in H^2(\Omega):\nabla u\cdot\nu=0 \mbox{ on } \Gamma \right\rbrace.
\end{aligned}
$$
By the spectral theory, there exists a complete orthonormal 
basis $\{\vp_k\}^{\infty}_{k=0}$ of $L^2(\Omega)$, where
$\vp_k$ satisfies the following eigenvalue problem 
\begin{equation}
\left\lbrace
\begin{aligned}
-\Delta \vp_k&=\lambda_k \, \vp_k,\quad\mbox{ in }\Omega,
\\[5pt]
\nabla \vp_k\cdot\nu&=0,\quad\quad\quad\mbox{ on }\Gamma.
\end{aligned}
\right.\label{Eq_Neumann}
\end{equation}
Therefore, we have that $\vp_k$ is the eigenfunction corresponding to eigenvalue $\lambda_k$
 for each $k\geq 0$, where one repeats each eigenvalue $\lambda_k$ according to its (finite) multiplicity:
$$
    0=\lambda_0< \lambda_1\leq\lambda_2\leq\lambda_3\leq \cdots \leq \lambda_k  \leq \cdots, \quad \text{$\lambda_k \rightarrow \infty$
     as $k \longrightarrow \infty$}.
$$
Moreover, it is not difficult to show that, 
$\vp_0=1/\sqrt{|\Omega|}$, that is a constant value, and 
\begin{equation}
\label{PHIMEDIA0}
   \int_\Omega \vp_k(x) \, dx= 0, \quad \text{for all $k\geq1$}.
\end{equation}
Then, we may write
$$
\begin{aligned}
   D(-\Delta_{N})&= \{u \in L^2(\Omega); \ \sum_{k= 1}^\infty \lambda_k^2 \, |\langle u, \vp_k \rangle|^2 < \infty \}, 
   \\[5pt]
   (-\Delta_{N})\, u&= \sum_{k= 1}^\infty \lambda_k \, \langle u, \vp_k \rangle \ \vp_k, 
   \quad \text{for each  $u \in D(-\Delta_{N})$}.
\end{aligned}
$$

\medskip
Now, applying the functional calculus, we define for each $s \in (0,1)$, the Neumann spectral fractional Laplacian operator, that is 
$$(-\Delta_{N})^{s}: D \big((-\Delta_{N})^{s}\big) \subset L^2(\Omega) \to L^2(\Omega),$$ 
given by
\begin{equation}
\label{DEFFRALAP}
\begin{aligned}
    (-\Delta_{N})^s u:&= \sum_{k= 1}^\infty \lambda^s_k \, \langle u, \vp_k \rangle \ \vp_k, 
\\
    D\big((-\Delta_{N})^s\big)&= \Big\{ u \in L^2(\Omega) : \; \sum^{\infty}_{k=1}
    \lambda^{2 s}_k \, \vert \langle u, \vp_k \rangle \vert^2<+\infty  \Big\}.
\end{aligned}
\end{equation}
Moreover, $D\big((-\Delta_{N})^s\big)$ is a Hilbert space, with the inner product
$$
      \langle u, v \rangle_{s}:=\langle u, v\rangle+\int_{\Omega} (-\Delta_{N})^{s} u(x) \; (-\Delta_{N})^{s} v(x) \, dx.
$$
In particular, the norm $\vert \cdot \vert_s$ is defined by
\begin{equation}
\label{eq:norm}
       \vert u \vert_s^2:= \Vert u \Vert^2_{L^2(\Omega)}+\Vert (-\Delta_{N})^su \Vert^2_{L^2(\Omega)}.
\end{equation}

Analogously, we define $\mathcal{K}_s: D\big(\mathcal{K}_s\big)= L^2(\Omega) \to L^2(\Omega)$ by
\begin{equation}
\label{DEFFRALAPINV}
\begin{aligned}
    \mathcal{K}_s u:&= \sum_{k= 1}^\infty \lambda^{-s}_k \, \langle u, \vp_k \rangle \ \vp_k, 
\\
    D\big(\mathcal{K}_s\big)&= \Big\{ u \in L^2(\Omega) : \; \sum^{\infty}_{k=1}
    \lambda^{-2s}_k \, \vert \langle u, \vp_k \rangle \vert^2<+\infty  \Big\}.
\end{aligned}
\end{equation}
The next proposition give us the main properties of the $(-\Delta_N)^{s}$, and $\mathcal{K}_s$ operators defined above.
\begin{proposition}
\label{THMCARAC}
Let $\Omega \subset \R^n$ be a bounded open set, $s \in (0,1)$, and consider $(-\Delta_N)^{s}$, and $\mathcal{K}_s$ the 
operators defined respectively by \eqref{DEFFRALAP} and \eqref{DEFFRALAPINV}. Then, we have:

\begin{enumerate}
\item[$(1)$]  $D(-\Delta_N) \subset D((-\Delta_N)^s)$, thus $D((-\Delta_N)^s)$ is dense in $L^2(\Omega)$.

\item[$(2)$] The operator   $(-\Delta_N)^{s}$ and $\mathcal{K}_s$ are self-adjoint.

\item[$(3)$] If $0 < s_{1}\leq s_{2} \leq1,$ then
$D\left((-\Delta_N)^{s_2}\right) \hookrightarrow D
\left((-\Delta_N)^{s_1}\right)$,
and 
$$D\left((-\Delta_N)^{s_2}\right) \mbox{ is dense in } D\left((-\Delta_N)^{s_1}\right).$$

\item[$(4)$] For any $\lambda> 0$, and $s \in [0,1]$, 
the operator $I_d + \lambda (-\Delta_N)^s$ is bijective from $D((-\Delta_N)^s)$ to $L^2(\Omega)$.
\end{enumerate}
\end{proposition} 
\begin{proof}
The proofs of items $(1)-(3)$ follow analogously to Proposition 2.1 in \cite{GHWN1}, hence we omit them. 
Let us show item (4), first we note that for any $u \in D((-\Delta_N)^s)$, we have 
\begin{equation}
  \langle (-\Delta_N)^{s} u, u \rangle= \sum_{k=1}^{\infty}\lambda_k^{s} \vert\langle u, \vp_k\rangle \vert^2\geq\lambda_1^{s}
  \sum_{k=1}^{\infty} \vert\langle u, \vp_k\rangle \vert^2\geq0,\label{eq:A-1}
\end{equation}
which implies for any $\lambda> 0$
\begin{equation}
 \label{acre}
    \Vert u + \lambda (-\Delta_N)^su \Vert_{L^2(\Omega)} \geq \Vert u \Vert_{L^2(\Omega)}.
\end{equation}
Therefore, the linear operator $I_d + \lambda(-\Delta_N)^s$ is injective. 
Moreover, for each $f \in L^2(\Omega)$ there exists $v \in D((-\Delta_N)^s)$, such that
$$
    v + \lambda (-\Delta_N)^s v= f.
$$
Indeed, it is enough to take
\begin{equation}
\label{INVERSA}
   v(x)= \sum^{\infty}_{k=0} \frac{\langle f, \vp_k \rangle}{1 + \lambda\lambda_k^s} \, \vp_k(x)
\end{equation}
and check that, $v \in D((-\Delta_N)^s)$ and satisfies the above equation. 
Therefore, $I_d + \lambda (-\Delta_N)^s$ is a bijective operator. 
\end{proof}

\begin{remark}\label{Re_inverse}
1. One remarks that, for each $\lambda> 0$ the operator 
$$
    \lambda \, I_d + (-\Delta_N)^s\,:\,D\big((-\Delta_N)^s\big)\to L^2(\Omega)
$$
is invertible. For the extremal case $(\lambda= 0)$, 
this assertion is false which is due to the fact that, $(-\Delta_N)^s$ is not injective in $D\big((-\Delta_N)^s\big)$.  

\medskip
2. Thanks to the above observation, we have that $(-\Delta_N)^s$ is not invertible in its domain. 
Then, the operator $\mathcal{K}_s$ could not be seen as the inverse of $(-\Delta_N)^s$. 
However, if we restrict the domain of the fractional Laplacian to a specific subset of $D\big((-\Delta_N)^s\big)$, 
we obtain the existence of the inverse $(-\Delta_N)^{-s}$.
% moreover its coincide with $\mathcal{K}_s$ in the restriction. 
\end{remark}
%%%%%%%%%%%%%%%%%%%%%%%%%%%%%%%%

Let us mention an important result, which help us to show the existence of solutions for the parabolic regularization of the system \eqref{eq.main}.
\begin{proposition}\label{Pro_inver_Lap_e_fracLap}
Given $v\in L^2(\Omega)$, then for all $s\in(0,1)$
$$
\,\mathcal{K}_s \big(\,I_d\,+\,(-\Delta_N)^{1-s}\,\big)^{-1}v\,\in\,D(-\Delta_N).
$$
\end{proposition}

\begin{proof} It follows from Proposition \ref{THMCARAC} together with the definition of the $\mathcal{K}_s$ operator.
\end{proof}
%%%%%%%%%%%%%%%%%%%%%%%%%%%%%%%%
\subsection{The inverse of the restricted NSFL operator}
%%%%%%%%%%%%%%%%%%%%%%%%%%%%%%%%

Here, we consider a subset of the domain $ D\big((-\Delta_N)^{s}\big)$, such that there exists the
inverse of the NSFL operator when restricted to this set. To this end, let us define for each $s \in [0,1)$
the following set
\begin{equation}
\label{2equalrepredominio}
\mathcal{H}^{2s}_N:=\left\lbrace u \in D\big((-\Delta_N)^{s}\big);  \int_{\Omega}u(x)\,dx=0\right\rbrace.
\end{equation}
Hence we have the following 
\begin{proposition}
\label{THMCARAC_1}
Under the conditions of Proposition \ref{THMCARAC}, it follows that:
%
%Let $\Omega \subset \R^n$ be a bounded open set, $s \in (0,1)$, and consider $(-\Delta_N)^{s}$, the 
%operators defined above. Then, we have:
\begin{enumerate}
\item[$(1)$] For all  $u \in \mathcal{H}^{2s}_N$, there exists $\alpha> 0$ such that
\begin{equation}
   \langle (-\Delta_N)^s u, u \rangle \geq \alpha^s \Vert u\Vert^2_{L^2(\Omega)}, \label{coercivity}
\end{equation}
where $\alpha$ is the coercivity constant of $(-\Delta_N)$. 

\item[$(2)$] The operator $
(-\Delta_N)^s$ is bijective from $\mathcal{H}^{2s}_N$ to $\mathcal{H}^0_N$. In particular, 
the inverse of the fractional Neumann spectral Laplacian, i.e. 
$((-\Delta_N)^{s})^{-1}$, exists.
\end{enumerate}
\end{proposition} 
\begin{proof}
1. First, for $u \in \mathcal{H}^{2s}_N$ we have
\begin{equation}
  \langle (-\Delta_N)^s u, u \rangle= \sum_{k=1}^{\infty}\lambda_k^{s} \vert\langle u, \vp_k\rangle \vert^2\geq\lambda_1^{s}
  \sum_{k=1}^{\infty} \vert\langle u, \vp_k\rangle \vert^2=\lambda_1^s\Vert u \Vert^2_{L^2(\Omega)}, \label{eq:A-1}
\end{equation}
where we have used in the last step, $\varphi_0= 1/ \sqrt{|\Omega|}$ and $\int_\Omega u(x) dx= 0$.

\medskip
2. From item (1), it follows that  $(-\Delta_N)^s$ is injective in $\mathcal{H}^{2s}_N$. Now, we show that $(-\Delta_N)^s$ is also surjective. Indeed, for any $u\in \mathcal{H}^{0}_N$ let $v$ be defined by
$$
 v:= \sum^{\infty}_{k=1} \lambda^{-s}_k \langle u,\vp_k\rangle \vp_k.
$$
Then, $v \in D\big((-\Delta_N)^s\big)$ and 
$$
\int_{\Omega}v(x)dx\,=\,\left\langle \sum^{\infty}_{k=1} \lambda^{-s}_k \langle u\, ,\,\vp_k\rangle \vp_k, 1\right\rangle= 0,
$$
where we have used that $\langle \vp_k,1\rangle=0 $, for each $k\geq1$. 
Consequently, $v \in \mathcal{H}^s_N$ and also $(-\Delta_N)^s v= u$. 
Then the operator$(-\Delta_N)^s$ is surjective, and thus $((-\Delta_N)^s)^{-1}$ exists. 
\end{proof}
%%%%%%%%%%
\begin{remark}
Applying Proposition \ref{THMCARAC_1}, it follows that for each $u\in \mathcal{H}^0_N$ the inverse 
of the NSFL is given by
\begin{equation}
\label{INVFNL}
((-\Delta_N)^s)^{-1}u =:(-\Delta_N)^{-s}u=\sum^{\infty}_{k=1} \lambda^{-s}_k \langle u,\vp_k\rangle \vp_k.
\end{equation}
From now own, we write $(-\Delta_N)^{-s}$ to denote the inverse of the NSFL operator,
whenever this makes sense. 
\end{remark}
%
%%%%%%%%%
\begin{proposition}
Let $(-\Delta_N)^{-s}$ be the 
operator defined by \eqref{INVFNL} for any fixed $s \in (0,1)$. Then, we have:
\begin{enumerate}
\item[$(1)$] The operator $(-\Delta_N)^{-s}$ is self-adjoint in $\mathcal{H}^{0}_N$.

\item[$(2)$] The operator $(-\Delta_N)^{-s}$ is continuous from $\mathcal{H}^{0}_N$ to itself.

\item[$(3)$] If $\sigma>0$ and $u \in \mathcal{H}^{2s}_N$, then $\mathcal{H}^{2(s+\sigma)}_N \ni v= (-\Delta_N)^{-\sigma}\,u$. 

%For each $s,\,\sigma>0$ and $u\in \mathcal{H}^{2s}_N$ we have$(-\Delta_N)^{-\sigma}\,u\in \mathcal{H}^{2s}_N$.
\end{enumerate}
\end{proposition}
\begin{proof}
The proof proceeds analogously to the Proposition 2.1 in \cite{GHWN1}, and hence
we omit it.
\end{proof}
%%%%%%%%%%%%%%%%%
%
%\begin{corollary}
%Let $v\in L^2(\Omega)$, we have
%$$
%\big(\,(-\Delta_N)\,+\,(-\Delta_N)^{s}\,\big)^{-1}(v-\left\langle v \right\rangle)\,\in\,D(-\Delta_N) ,
%$$
%for all $s\in(0,1)$.
%\end{corollary}
%%%%%%%%%%%%%%%%%%%
\subsection{Some auxiliary results}
%%%%%%%%%%%%%%%%%%%

First, we recall that using the language of semigroups, as introduced in \cite{Stinga} (see also \cite{ST}), 
one can check that $\left(-\Delta_{N}\right)^{s}$ is indeed a nonlocal operator. In fact, the NSFL is also given by
$$
\left(- \Delta_{N}\right)^{s} u(x)=\frac{1}{\Gamma(-s)} \int_{0}^{\infty}\left(e^{t \Delta_{N}} u(x)-u(x)\right) \frac{d t}{t^{1+s}}, \quad x \in \Omega,
$$
where $e^{t \Delta_{N}} u(x)$ is the heat diffusion semigroup generated by the Neumann Laplacian acting on $u$. 

\medskip
Now, the aim  is to characterize  the space $D((-\Delta_N)^{s})$. To begin, we study $D((-\Delta_N)^{1/2})$, indeed, by using the 
$L^{2}$ normalization and the weak formulation of the equation \eqref{Eq_Neumann}, we see that $\left\|\varphi_{k}\right\|_{H^{1}(\Omega)}^{2}=1+\lambda_{k} .$ 
It is easy to check that $\left\{\varphi_{k}\right\}_{k \in \mathbb{N}_{0}}$ is also an orthogonal basis of $H^{1}(\Omega)$.
Hence, we find
$$
H^{1}(\Omega)=\Big\{u \in L^{2}(\Omega):\|u\|_{H^{1}(\Omega)}^{2}=\sum_{k=0}^{\infty}\left(1+\lambda_{k}\right)\left|\left\langle u, \varphi_{k}\right\rangle\right|^{2}<\infty\Big\}.
$$

Therefore $H^1(\Omega)=D((-\Delta_N)^{1/2})$. In particular, from \eqref{eq:norm} we have
\begin{equation}
\int_{\Omega} (-\Delta_{N})^{1/2} u(x) \ (-\Delta_{N})^{1/2} u(x) \ dx=\int_{\Omega}\nabla u(x) \cdot\nabla u(x) \ dx\label{normequi}
\end{equation}
for all $u\in H^1(\Omega)$. Similarly, we have the following
\begin{proposition}
\label{represendoDominio}
%\label{eq:dominequiv}
Let $\Omega \subset \Bbb{R}^n$ be a bounded open set. Then for any $s \in (0,1)$, 
\begin{equation}
\mathcal{H}^{s}_N=\left\lbrace u \in H^s(\Omega):  \int_{\Omega}u(x)\,dx=0\right\rbrace.
\label{equalrepredominio}
\end{equation}
In particular, for each $u\in \mathcal{H}^{s}_N$, there exist $C_1, C_2> 0$ such that
$$
\begin{aligned}
   C_1 \iint_{\Omega\times\Omega}\frac{|u(x)-u(y)|^2}{|x-y|^{n+2s}}\,dx\,dy
  & \leq \|\,(-\Delta_N)^{s/2}u\,\|^2_{L^2(\Omega)} 
  \\
  &\leq C_2 \iint_{\Omega\times\Omega}\frac{|u(x)-u(y)|^2}{|x-y|^{n+2s}}\,dx\,dy.
\end{aligned}
$$
\label{eq:dominequiv} 
\end{proposition}
\begin{proof}
See Theorem 2.5 in \cite{SV}, and Lemma 7.1 in \cite{CS}.  
\end{proof}

Here and subsequently, we denote for each $s \in (0,1)$ the operator
$\mathcal{H}_s= \clg{K}_s^{1/2}$. 
Then, we consider the following 
\begin{lemma}\label{lem:Kboundary}
Let $\;\Omega \subset \R^n$ be a bounded open set, $s\in(0,1)$ and $ u \in D((-\Delta_N)^{1-s})$, then
$\mathcal{K}_su \in D(\,-\Delta_N\,)$. In particular, we have  in trace sense
$$
\nabla \mathcal{K}_su \cdot \nu =0 \mbox{ on } \Gamma.
$$
\end{lemma}
\begin{proof}
The proof follows directly from Proposition \ref{THMCARAC}, item (3).
\end{proof}

\begin{proposition}
\label{proKuHu}
Let $\;\Omega \subset \R^n$ be a bounded open set and $0<s<1$.
\begin{enumerate}
\item[$(1)$] 
There exists a constant $C_{\Omega}> 0$, such that, for each $u \in H^1(\Omega)$
\begin{equation}
\label{nablaKunablau}
   \int_{\Omega}|\nabla\mathcal{K}_su(x)|^2 \ dx
   \leq C_{\Omega}\int_{\Omega}|\nabla u(x)|^2 \ dx.
\end{equation}
Similarly, for each $u\in H^1(\Omega)$, then $\nabla\mathcal{H}_s u\in \mathbf{L}^2(\Omega)$ and
\begin{equation}
\label{nablaHunablau}
   \int_{\Omega}|\nabla\mathcal{H}_s u(x)|^2 \ dx
   \leq C_{\Omega}^{1/2} \int_{\Omega}|\nabla u(x)|^2 \ dx.
\end{equation}

\item[$(2)$] If $u\in H^1(\Omega)$, then
\begin{equation}
\label{pro:gradKugradu}
    \int_{\Omega}\nabla \mathcal{K}_s u(x)\cdot \nabla u(x)dx = \int_{\Omega}\vert \nabla \mathcal{H}_s u(x) \vert^2 \ dx.
\end{equation}
\end{enumerate}
\end{proposition}
\begin{proof}
1. To show item 1, we use the equivalence between
 $D((-\Delta_N)^{1/2})$ and $H^1(\Omega)$ (see \eqref{normequi}).
 Then, we have
$$
\begin{aligned}
\int_{\Omega}|\nabla\mathcal{K}_su(x)|^2dx&=\sum_{k=1}^{\infty}\lambda_k\vert \langle \mathcal{K}_s u, \vp_k \rangle \vert^2
=\sum_{k=1}^{\infty}\lambda_k\vert\lambda_k^{-s}\langle u, \vp_k \rangle\vert^2
\\[5pt]
&\leq \lambda_1^{-2s}\sum_{k=1}^{\infty}\lambda_k\vert\langle u, \vp_k \rangle\vert^2=\lambda_1^{-2s}\int_{\Omega}|\nabla u(x)|^2dx< \infty,
\end{aligned}
$$ 
and analogously for $\nabla \clg{H}_s u$.

\medskip
2. Now, we show item 2. Since $u\in H^1(\Omega)$, it is enough to consider that $u \in D(\,-\Delta_N\,)$
and thus apply a standard density argument. First, we integrate by parts to obtain
$$
    \int_{\Omega}\nabla \mathcal{K}_su(x)\cdot\nabla u(x)\,dx
   =\int_{\Omega} (-\Delta_N) \mathcal{K}_su(x)u(x)\,dx
   =\int_{\Omega}(-\Delta_N)^{1-s}u(x)u(x)\,dx,
$$
where we have used the definition of $\mathcal{K}_su$ and $\nabla\mathcal{K}_su\cdot\nu=0$ on $\Gamma$.
Due to the fractional Laplacian being self-adjoint, it follows that
$$
\int_{\Omega}\nabla \mathcal{K}_su(x)\cdot\nabla u(x)\,dx=\int_{\Omega}|(-\Delta_N)^{(1-s)/2}u(x)|^2\,dx.
$$
Therefore, using the equivalence norm  \eqref{normequi} together with the definition of $\mathcal{H}_s u$, we have
$$
\int_{\Omega}\nabla \mathcal{K}_su(x)\cdot\nabla u(x)\,dx=\int_{\Omega}|\nabla \mathcal{H}_s u|^2\, dx.
$$
\end{proof}

%%%%%%%%%%%%%%%%%%
\section{On a Perturbed Problem}
\label{ParabolicApp}
%%%%%%%%%%%%%%%%%%
The aim of this section is to introduce and study the properties of a perturbed system 
associated to \eqref{eq.main}. More precisely, given $\varepsilon>0$ we consider the 
following fractional parabolic-elliptic system
\begin{equation}
\label{eq_regulation}
\begin{cases}
\partial_{t} u_{\varepsilon} + \operatorname{div} 
\left(g(u_{\varepsilon})\nabla\mathcal{K}_s c_{\varepsilon} \right)= \varepsilon \Delta u_{\varepsilon}, 
\quad &\text{in $(0,\infty) \times \Omega$},  
\\[5pt]
(-\Delta_N)^{1-s} c_{\varepsilon}+c_{\varepsilon}= u_{\varepsilon}, \quad & \text{in $\Omega$}, 
\\[5pt]
u_{\varepsilon}(0)= u_{0,\varepsilon}, \hspace{25pt}  &\text{in $\Omega$},
\\[5pt]
\nabla \mathcal{K}_sc_{\varepsilon} \cdot\nu= 0 \text{ and } \nabla u_{\varepsilon} \cdot\nu= 0, \hspace{30pt} &\text{on $\Gamma$},
\end{cases}
\end{equation}
where $u_{0,\varepsilon}$ is a regularized initial data, 
such that 
\begin{equation}
\label{U0}
\text{$u_{0,\varepsilon} \to u_0 $ strongly in $L^{1}(\Omega)$
as $\varepsilon \to 0$,
and $\|u_{0,\varepsilon}\|_{L^\infty} \leq \|u_{0}\|_{L^\infty}$.}
\end{equation}
Then, we show existence and uniqueness of $(u_\varepsilon,c_\varepsilon)$. To this end, 
we are going to apply the Banach Fixed Point Theorem to prove the local in time existence of solution to \eqref{eq_regulation},
and thus applying a contradiction argument we extend it to be global in time. Since \eqref{eq_regulation} is a fractional non-standard 
parabolic-elliptic system we present the proof in details. 
To begin, we consider the following
\begin{lemma} 
\label{LEMMA31}
Let $s \in (0,1)$ be fixed and $\tilde{u}\in L^{\infty}((0,\infty) \times \Omega)$.
Then for each $u_{0,\varepsilon} \in L^2(\Omega)$ the fractional parabolic-elliptic system
\begin{equation}
\label{Frac_para_ellip}
\begin{cases}  
\partial_{t} u_{\varepsilon}-\varepsilon \Delta u_{\varepsilon}=- \dive 
\left(g(\tilde{u})\nabla\mathcal{K}_s c_{\varepsilon} \right) , 
\quad& \text{in $(0,\infty) \times \Omega$},  
\\[5pt]
(-\Delta_N)^{1-s} c_{\varepsilon}+c_{\varepsilon}= \tilde{u}, \quad& \text{in $\Omega$}, 
\\[5pt]
u_{\varepsilon}(0)= u_{0,\varepsilon}, \hspace{25pt}  &\text{in $\Omega$},
\\[5pt]
\nabla \mathcal{K}_sc_{\varepsilon} \cdot \nu= 0, \text{ and } \nabla u_{\varepsilon} \cdot \nu= 
0 \hspace{30pt} &\text{on $\Gamma$},
\end{cases}
\end{equation}
admits a unique weak solution 
$$
    (u_\varepsilon, c_\varepsilon) \in L^2\big((0,\infty);H^1(\Omega)\big)\cap C\big([0,\infty);L^2(\Omega)\big) \times L^{\infty}\big((0,\infty); D\big((-\Delta_N)^{1-s}\big)\big).
$$    
\end{lemma}
\begin{proof}
1. Given $\tilde{u} \in L^{\infty}((0,\infty) \times \Omega)$ from \eqref{INVERSA} and Remark \ref{Re_inverse}, it is possible to write, for almost all $t> 0$, the chemoattractant density as follows
\begin{equation}
\begin{aligned}
c_{\varepsilon}(t)= \big(\,I_d\,+\,(-\Delta_N)^{1-s}\,\big)^{-1} \tilde{u}(t)=
\,\sum_{k=0}^{\infty}\,\frac{1}{1+\lambda_k^{1-s}}\,\langle \tilde{u}(t),\varphi_k\rangle\,\varphi_k
\end{aligned}\label{Eq_c_epsilon}
\end{equation}
In particular $c_\varepsilon\in L^{\infty}\big((0,\infty); D\big((-\Delta_N)^{1-s}\big)\big)$ and, for almost all $t>0$ and 
almost everywhere $x \in \Omega$, satisfies the equation
$$
(-\Delta_N)^{1-s}c_\varepsilon(t,x) + c_\varepsilon(t,x)= \tilde{u}(t,x).
% \quad \text{ in } L^2(\Omega).
$$ 
Moreover, due to Proposition \ref{Pro_inver_Lap_e_fracLap} we have that $\mathcal{K}_sc_\varepsilon(t)\in D(-\Delta_N)$ 
for a.a. $t> 0$. In particular, $\nabla\mathcal{K}_sc_\varepsilon\cdot\nu=0$ on $\Gamma$, which is due to the definition of $D(-\Delta_N)$. 
Consequently, the second equation in \eqref{Frac_para_ellip} is solved. The uniqueness follows easily from the linearity of the equation. 

\medskip
2. Now, we show the existence of solution for the first equation in \eqref{Frac_para_ellip}. First, since we have the explicit form of $c_\varepsilon$, we may write
\begin{eqnarray}
\mathcal{K}_s c_\varepsilon&=&\mathcal{K}_s \big(\,I_d\,+\,(-\Delta_N)^{1-s}\,\big)^{-1}\tilde{u}\label{eq_oper_1}
\\[5pt]
&=&\sum_{k=1}^{\infty}\frac{\lambda_k^{-s}}{1+\lambda_k^{1-s}}\left\langle\tilde{u},\varphi_k\right\rangle\,\varphi_k
\,\,=\,\,\sum_{k=1}^{\infty}\frac{1}{\lambda_k+\lambda_k^{s}}\left\langle\tilde{u},\varphi_k\right\rangle\,\varphi_k,\label{eq_sum_1}
\end{eqnarray}
where \eqref{eq_sum_1} follows by the definition of $\mathcal{K}_s$ together with \eqref{Eq_c_epsilon}. 
Observe that, $\mathcal{K}_sc_\varepsilon$ is an operator which depends on $\tilde{u}$, so we define
$$
\mathcal{L}\tilde{u}\,:=\,\mathcal{K}_sc_\varepsilon\,=\,\sum_{k=1}^{\infty}\frac{1}{\lambda_k+\lambda_k^{s}}\left\langle\tilde{u},\varphi_k\right\rangle\,\varphi_k,
$$
where the last assertion is obtained by \eqref{eq_oper_1}.  Then the first equation in \eqref{Frac_para_ellip} together with the 
initial-boundary condition can be written as follows
\begin{equation}
\label{ecua153}
\left\{
\begin{aligned}
 \partial_{t} u_\varepsilon - \varepsilon \Delta u_\varepsilon
  & = - \dive\big(g(\tilde{u})\ \nabla\mathcal{L}\tilde{u}\big), \quad \text{in $(0,\infty) \times \Omega$}, 
 \\[5pt]
 u_\varepsilon(0) &= u_{0,\varepsilon}, \hspace{20pt} \text{in $\Omega$},
  \\[5pt]
\nabla u_\varepsilon \cdot \nu  & =0, \hspace{32pt} \text{on $\Gamma$}.
\end{aligned}
\right.
\end{equation}
Therefore, it is enough to show the existence of a unique solution for \eqref{ecua153}. 

We claim that, for any $\tilde{u} \in L^{\infty}((0,\infty) \times \Omega)$, $$\dive\big(g(\tilde{u})\ \nabla\mathcal{L}\tilde{u}\big)\in L^{\infty}((0,\infty);H^{-1}(\Omega)).$$ 
Indeed, from \eqref{eq_sum_1}, the definition of $\mathcal{L}\tilde{u}$ and the equivalence norm (see \eqref{normequi}), we have for almost all $t> 0$
$$
\begin{aligned}
\int_{\Omega} |\nabla\mathcal{L}(\tilde{u}(t,x))|^2 dx&= \int_{\Omega}|(-\Delta_N)^{1/2}\mathcal{L}(\tilde{u}(t,x))|^2  dx
\\[5pt]
&=\,\sum_{k=1}^{\infty}\lambda_k\left(\frac{1}{\lambda_k+\lambda_k^{s}}\right)^2|\left\langle\tilde{u}(t),\varphi_k\right\rangle|^2
\\[5pt]
&\,=\,\sum_{k=1}^{\infty}\lambda_k^{-1}\left(\frac{\lambda_k}{\lambda_k+\lambda_k^{s}}\right)^2|\left\langle\tilde{u}(t),\varphi_k\right\rangle|^2
\\[5pt]
&\,\leq\,\lambda_1^{-1}\sum_{k=1}^{\infty}|\left\langle\tilde{u}(t),\varphi_k\right\rangle|^2\,\leq\,\lambda_1^{-1} \int_\Omega |\tilde{u}(t,x)|^2 dx
\\[5pt]
&\,\leq\,\lambda_1^{-1}\,|\Omega|\,\|\tilde{u}\|^2_{L^{\infty}},
\end{aligned}
$$
where $|\Omega|$ is the Lebesgue measure of $\Omega$ and $\lambda_1$ is the first eigenvalue of \eqref{Eq_Neumann}. 
Therefore, we obtain $g(\tilde{u})\nabla\mathcal{L}\tilde{u}\in L^{\infty}\big((0,\infty);\mathbf{L}^2(\Omega)\big)$, 
from which follows the claim.

Finally, applying a standard parabolic theory (see Lions, Magenes \cite{LionsMagenes2}),
there exists a unique weak solution $u_\varepsilon \in L^2\big((0,\infty);H^1(\Omega)\big)\cap C\big([0,\infty);L^2(\Omega)\big)$ 
of \eqref{ecua153},
and thus the proof is complete. 
\end{proof}

From the proof of the previous lemma, one observes that the fractional parabolic-elliptic system \eqref{eq_regulation} can be decoupled. Therefore, 
equivalently to show existence and uniqueness of solution $(u_\varepsilon,c_\varepsilon)$ for \eqref{eq_regulation}, we study 
the following initial-boundary Neumann value problem
\begin{equation}
\label{Eq_PR_modificado}
\left\{
\begin{aligned}
\partial_{t} u_{\varepsilon}-\varepsilon \Delta u_{\varepsilon}&= -\dive\big(\,g(u_{\varepsilon})\nabla\,\mathcal{L}(u_{\varepsilon}) \,\big), \quad \text {in $\Omega_{T}$}, 
\\[5pt]
\nabla u_\varepsilon\cdot\nu&= 0, \quad \quad  \text{on $\Gamma$},
\\[5pt]
u_\varepsilon(0)&=u_{0,\varepsilon}, \quad \text {in $\Omega$}. 
\end{aligned}
\right.
\end{equation}

%%%%%%%%%%%%%%%%%%%%%%%%%%%
\subsection{Local existence of solution} 
%%%%%%%%%%%%%%%%%%%%%%%%%%%

For convenience, let us denote for any $T> 0$,
$$W(T)= L^2\big((0,T);H^1(\Omega)\big)\cap C\big([0,T);L^2(\Omega)\big),$$
and also $\Omega_T= (0,T) \times \Omega$. Then, we consider the following 

\begin{theorem}[Local existence]
\label{th_local}
Given $u_{0,\epsilon}$ satisfying \eqref{U0}, then there exists a positive time $T=T(u_{0})$, such that 
the problem \eqref{Eq_PR_modificado} admits a unique weak
solution 
\[
u_\varepsilon \in L^{\infty}\left(\Omega_{T}\right) \cap W(T).
\]
\end{theorem}
\begin{proof}
1. Hereupon,  we denote by $B_R^T$ the following set 
$$
B_R^T:=\left\lbrace \tilde{u}\in L^{\infty}(\Omega_T): \|\tilde{u}\|_{L^{\infty}(\Omega_T)}\leq R\right\rbrace,
$$
where $T>0$ and $R>0$ are chosen a posteriori.
Also we define the mapping $\mathcal{T}:B_R^T \rightarrow W(T)$, $\tilde{u} \mapsto u_\varepsilon= \mathcal{T}(\tilde{u})$ where 
$u_\varepsilon$ is the unique solution of \eqref{ecua153} (for each $\varepsilon>0$ fixed), that is to say, for each 
$\tilde{u} \in L^{\infty}(\Omega_{T})$ and $t \in [0,T)$, we may write
\begin{equation}
\label{eq_math_T}
\begin{aligned}
u_\varepsilon(t,x) \equiv \mathcal{T}(\tilde{u})(t,x)&=  \int_\Omega K(t,x,y) u_{0,\varepsilon}(y) \ dy  
\\[5pt]
&+\int_{0}^{t} \int_{\Omega} g(\tilde{u}(\tau,y)) \,\nabla_{\!\!y} K(t-\tau,x,y) \cdot \nabla\mathcal{L}\tilde{u}(\tau,y)\,dy\, d \tau, 
\end{aligned}
\end{equation}
where $K(t,x,y)$ is the heat kernel of $(-\Delta)$ with Neumann boundary condition, namely (see  \cite{DAVIES}, Theorem 2.1.4) 
$$
K(t,x,y)=\sum_{k=0}^{\infty}e^{-t\varepsilon\lambda_k}\varphi_k(x)\,\varphi_k(y).
$$

\medskip
2.  Now, we show that $\mathcal{T}(\tilde{u})\in B_R^T$ for each $\tilde{u}\in B_R^T$. Indeed, from  \eqref{eq_math_T} and 
applying H\"older inequality, we obtain
\begin{equation}
\begin{aligned}
|u_\varepsilon(t,x)| &\leq C \|u_{0}\|_{L^{\infty}(\Omega)}
\\[5pt]
&+\int_0^t\|\nabla K(t-\tau,x,\cdot)\|_{L^1(\Omega)}\,\|g(\tilde{u}(\tau))\|_{L^\infty(\Omega)}\,\|\nabla\mathcal{L}\tilde{u}(\tau)\|_{L^\infty(\Omega)} d\tau
\\[5pt]
\,&\leq C \|u_{0}\|_{L^{\infty}(\Omega)} + C_1 \int_0^t(t-\tau)^{-1/2}\|g(\tilde{u}(\tau))\|_{L^\infty(\Omega)}\,\|\nabla\mathcal{L}\tilde{u}(\tau)\|_{L^\infty(\Omega)}\,d\tau,
\end{aligned}\label{eq_norm_inft_T}
\end{equation}
where $C,C_1$ are positive constants and we have used 
estimates of heat kernels with Neumann
boundary conditions, (see \cite{YangZhang}, Lemma 3.3). 

To follow, from \eqref{eq_sum_1} let $G_s(x,y)$ be the kernel of $\mathcal{L}$, that is to say, for each
$s \in (0,1)$, we define
$$
G_s(x,y):= \sum_{k=1}^{\infty}\frac{1}{\lambda_k+\lambda_k^s}\,\varphi_k(x)\,\varphi_k(y), \quad (x,y \in \Omega). 
$$
Consequently, $\nabla\mathcal{L}\tilde{u}$ could be written as follows
$$
\nabla\mathcal{L}\tilde{u}(t,x)= \int_{\Omega}\nabla G_s(x,y)\,\tilde{u}(t,y)\,dy,
$$
for almost everywhere $t \in (0,T)$, $x \in \Omega$.  
Therefore, applying Theorem 6.4 in  \cite{ALP} (see also \cite{CS}, Proposition 5.2),
we have
$$
\|\nabla\mathcal{L}\tilde{u}(t)\|_{L^\infty(\Omega)} \leq C_2 \  \|\tilde{u}(t)\|_{L^\infty(\Omega)},
$$
where $C_2$ is a positive constant. 
From the above estimate and \eqref{eq_norm_inft_T}, we obtain
$$
\begin{aligned}
|u_\varepsilon(t,x)| &\leq C \|u_0\|_{L^{\infty}(\Omega)}
\\[5pt]
&+ C_1 C_2 \int_0^t(t-\tau)^{-1/2}\|g(\tilde{u}(\tau))\|_{L^\infty(\Omega)}\,\|\tilde{u}(\tau)\|_{L^\infty(\Omega)}\,d\tau
%\\[5pt]
%\,&\leq\,c\|u_{0,\varepsilon}\|_{L^{\infty}(\Omega)}+c_2\int_0^t(t-\tau)^{-1/2}\left(\|\tilde{u}(\tau)\|_{L^\infty(\Omega)}+\|\tilde{u}(\tau)\|_{L^\infty(\Omega)}^2\right)\,\|\tilde{u}(\tau)\|_{L^\infty(\Omega)}\,d\tau
\\[5pt]
\,&\leq\, C \|u_0\|_{L^{\infty}(\Omega)} + C_1 C_2 \, R(R+R^2) \ 2 T^{1/2}.
\end{aligned}
$$
Therefore, taking $R= \dfrac{6}{5} C \|u_0\|_{L^{\infty}(\Omega)}$ and $T= \dfrac{1}{(12 C_1 C_2(R+R^2))^2}$,
it follows that  $\mathcal{T}(\tilde{u}) \in B_R^T$.

\medskip
3. Finally, we prove that $\mathcal{T}$ is a contraction on $B_R^T$. Indeed, we consider $\tilde{u}_i\in B_R^T$, ($i=1,2$), and similarly to item (2) we have
$$
\begin{aligned}
|&\mathcal{T}(\tilde{u}_1)(t,x) - \mathcal{T}(\tilde{u}_2)(t,x)|
\\[5pt]
&\,\leq\,\int_{0}^t\|\nabla K(t-\tau,x,\cdot)\|_{L^1(\Omega)} \|g(\tilde{u}_1)\nabla\mathcal{L}(\tilde{u}_1)(\tau)- g(\tilde{u}_2) \nabla\mathcal{L}(\tilde{u}_2)(\tau)\|_{L^{\infty}(\Omega)} d\tau
\\[5pt]
%&\,\leq  C_1\int_{0}^t(t-\tau)^{-1/2}\|\tilde{u}_1(1-\tilde{u}_1)\nabla\mathcal{L}(\tilde{u}_1)(\tau)-\tilde{u}_2(1-\tilde{u}_2)\nabla\mathcal{L}(\tilde{u}_2)(\tau)\|_{L^{\infty}(\Omega)}\,d\tau
%\\[5pt]
%&\,\leq\,c_1\int_{0}^t(t-\tau)^{-1/2}\|\tilde{u}_1(1-\tilde{u}_1)\nabla\mathcal{L}(\tilde{u}_1)(\tau)-\tilde{u}_2(1-\tilde{u}_1)\nabla\mathcal{L}(\tilde{u}_1)(\tau)\|_{L^{\infty}(\Omega)}\,d\tau
%\\[5pt]
%&\,+\,c_1\int_{0}^t(t-\tau)^{-1/2}\|\tilde{u}_2(1-\tilde{u}_1)\nabla\mathcal{L}(\tilde{u}_1)(\tau)-\tilde{u}_2(1-\tilde{u}_2)\nabla\mathcal{L}(\tilde{u}_1)(\tau)\|_{L^{\infty}(\Omega)}\,d\tau
%\\[5pt]
%&\,+\,c_1\int_{0}^t(t-\tau)^{-1/2}\|\tilde{u}_2(1-\tilde{u}_2)\nabla\mathcal{L}(\tilde{u}_1)(\tau)-\tilde{u}_2(1-\tilde{u}_2)\nabla\mathcal{L}(\tilde{u}_2)(\tau)\|_{L^{\infty}(\Omega)}\,d\tau 
%\\[5pt]
&\,\leq\,3 C_1 C_2 \ (1+R) R \,T^{1/2}\|\tilde{u}_1-\tilde{u}_2\|_{L^{\infty}(\Omega_T)} \,\leq\,\frac{1}{4}\|\tilde{u}_1-\tilde{u}_2\|_{L^{\infty}(\Omega_T)}.
\end{aligned}
$$
Therefore, the mapping $\mathcal{T}$ is a contraction, and hence
 we can apply the Banach Fixed Point Theorem. Thus
 $\mathcal{T}$ has a fixed point, which is by construction the unique solution of  
 \eqref{Eq_PR_modificado}. 
\end{proof}
\begin{remark}
\label{rem_TM}
Let $T_M$ be the maximal time of existence of the solution $u_\varepsilon$ for the problem \eqref{Eq_PR_modificado}. 
If $T_M<\infty$, then there exists an increase sequence $\{t_j\}_{j=1}^\infty$, such that,
$t_j\to T_M^{-}$ as $j\to\infty$ and
$$
\lim_{j\to\infty}\|u_\varepsilon(t_j,\cdot)\|_{L^{\infty}(\Omega)}=+\infty.
$$
\end{remark}

%%%%%%%%%%%%%%%%%%%%%%
\subsection{Global existence of solution}
%%%%%%%%%%%%%%%%%%%%%%

The main issue of this section is to show that, under the condition 
\eqref{ID} for the initial data $u_0$, we obtain (by contradiction) global in time existence of
solution of the problem \eqref{eq_regulation}. To this end, we show first the uniform boundedness 
of $(u_\varepsilon, c_\varepsilon)$.  
\begin{proposition}
\label{pro_bound_u_c}
Let $u_0 \in L^{\infty}(\Omega)$ be satisfying \eqref{ID}, and consider
$$
    (u_\varepsilon, c_\varepsilon) \in W(T_M) \cap L^{\infty}(\Omega_{T_M}) \times L^{\infty}\big((0,T_M); D\big((-\Delta_N)^{1-s}\big)\big)
$$    
the unique weak solution of \eqref{eq_regulation}. Then, it satisfies
\begin{eqnarray}
0 \leq u_{\varepsilon}(t, x) \leq 1,\quad  \text { a.e. in }[0, T_M) \times \Omega \label{Pro_eq_bound_u},
\\[5pt]
0 \leq c_{\varepsilon}(t, x) \leq 1,\quad  \text { a.e. in }[0, T_M) \times \Omega \label{Pro_eq_bound_c}. 
\end{eqnarray}
\end{proposition}
\begin{proof}
1. First, we observe that,
$\dive\left(g(u_\varepsilon)\nabla\mathcal{L}u_\varepsilon \right)\,\in\, L^2\left((0,T_M);L^2(\Omega)\right)$. 
Therefore, from equation \eqref{Eq_PR_modificado} and standard parabolic regularity theory, we obtain
$$
u_{\varepsilon}\in C\left([0,T_M);H^1(\Omega)\right)\cap L^2\left((0,T_M);H^2(\Omega)\right)
\,\,\text{ and }\,\,\,
\partial_t u_\varepsilon \in L^2\left(\Omega_{T_M}\right).
$$
Consequently, the pair $(u_\varepsilon, c_\varepsilon)$ satisfies the partial differential equations in \eqref{eq_regulation}
in the strong sense, that is, for almost all $(t,x) \in \Omega_{T_M}$. 

\medskip
2. To show \eqref{Pro_eq_bound_u}, we multiply the first equation in \eqref{eq_regulation} by $\varphi_{\delta}^{\prime}(u_\varepsilon),$ where
for $\delta> 0$
\[
\varphi_{\delta}(z)=\left\{\begin{array}{ll}
\left((z-1)^{2}+\delta^{2}\right)^{1 / 2}-\delta, & \text { for } z \geq 1, 
\\[5pt]
\quad 0 ,& \text { for } z \leq 1.
\end{array}\right. 
\]
Then, we obtain for each $t \in (0,T_M)$ 
\[
\begin{aligned}
\int_{\Omega} \varphi_{\delta}(u_\varepsilon(t,x)) d x
&-\iint_{\Omega_t}\varphi_{\delta}^{\prime \prime}(u_\varepsilon(\tau,x)) \ g(u_\varepsilon(\tau,x)) \nabla\mathcal{K}_sc_\varepsilon(\tau,x)  \cdot \nabla u_\varepsilon(\tau,x) \  d x d \tau 
\\[5pt]
&+\varepsilon \iint_{\Omega_t}|\nabla u_\varepsilon(\tau,x)|^{2} \varphi_{\delta}^{\prime \prime}(u_\varepsilon(\tau,x)) d x d \tau= 0,
\end{aligned}
\]
where we have used that $0\leq u_0\leq1$, and the boundary condition in \eqref{eq_regulation}. On the other hand, one observes that 
\[
\begin{aligned}
-\varphi_{\delta}^{\prime \prime}(u_\varepsilon)g(u_\varepsilon)\nabla\mathcal{K}_sc_\varepsilon\cdot \nabla u_\varepsilon&+\varepsilon|\nabla u_\varepsilon|^{2} \varphi_{\delta}^{\prime \prime}(u_\varepsilon)
 \\[5pt]
& \geq\left\{-|g(u_\varepsilon)||\nabla\mathcal{K}_sc_\varepsilon||\nabla u_\varepsilon|+\varepsilon|\nabla u_\varepsilon|^{2}\right\} \varphi_{\delta}^{\prime \prime}(u_\varepsilon) 
\\[5pt]
& \geq-\frac{1}{4 \varepsilon}u_\varepsilon^2(u_\varepsilon-1)^{2}|\nabla\mathcal{K}_sc_\varepsilon|^2 \varphi_{\delta}^{\prime \prime}(u_\varepsilon)
\\[5pt]
& \geq-\frac{\delta}{4 \varepsilon}u_\varepsilon^2|\nabla\mathcal{K}_sc_\varepsilon|^2,
\end{aligned}
\]
where we have used $g(u_\varepsilon)=u_\varepsilon(1-u_\varepsilon)$ and $(u_\varepsilon-1)^2\varphi_{\delta}^{\prime \prime}(u_\varepsilon)\leq\delta$. 
Consequently, 
\[
\int_{\Omega} \varphi_{\delta}(u_\varepsilon(t)) d x\leq\frac{\delta}{4 \varepsilon}\int_{\Omega_t}u_\varepsilon^2\,|\nabla\mathcal{K}_sc_\varepsilon|^2\,dx\,d\tau,
\]
and passing to the limit as $\delta\to0^+$, we have
$$
\int_\Omega|u_\varepsilon(t,x)-1|^{+} \ dx\leq 0.
$$
Thus for a.e. $(t,x) \in \Omega_{T_M}$, $|u_\varepsilon(t,x)-1|^{+}= 0$ and similarly we show that $|u_\varepsilon(t,x)|^{-}= 0$, 
from which follows \eqref{Pro_eq_bound_u}.

%\begin{comment}
%\medskip
%3. It remains to show \eqref{Pro_eq_bound_c}. We claim that $c_\varepsilon\leq 1$ in $[0,T_M)\times\Omega$. On contrary, suppose that
%$$
%A=\left\lbrace\,(t,x)\in\Omega_{T_M};\,1<c_\varepsilon(t,x)\,\right\rbrace
%$$
%has positive measure. Then integration the second equation in \eqref{eq_regulation} over $A$, we obtain
%$$
%\int_A c_\varepsilon(t,x)dtdx=\int_A u_\varepsilon(t,x)dtdx.
%$$
%since $\partial A=\left\lbrace (t,x);\, c_{\varepsilon}(t,x)=1\right\rbrace$ together with $\mathcal{K}_sc=0$ on $\partial A$. On the other hand, due to the fact that $0\leq u_\varepsilon\leq 1$, we have
%$$
%\int_A c_\varepsilon(t,x)dtdx\leq |A|,
%$$
%which is a contradiction. Therefore, for a.e. $(t,x) \in \Omega_{T_M}$, $c_\varepsilon(t,x)\leq 1$. By an analogous argument, we obtain that $0\leq c_\varepsilon(t,x)$, which finish the proof.
%\end{comment}

\medskip
3. Now, we recall that $\Omega\subset\R^n$ (with $n=1,\,2,\,3$), then by the regularity of $c_\varepsilon$ together with Morrey's Inequality, 
we have that $c_\varepsilon\in L^{\infty}\big( 0,T_M;\, C(\bar{\Omega})\big)$, thus the $\inf c_\varepsilon$
is finite, where the infimum is taking over $(0,T_M)\times\bar{\Omega}$.

\medskip
On the other hand, without loss of generality, let $(t_0,x_0)\in(0,T_M)\times\Omega$
 be a point where $c_\varepsilon(\cdot,\cdot)$ is a minimum, which is to say
$$
   c_\varepsilon(t,x)\geq c_\varepsilon(t_0,x_0) \quad \text{for each $(t,x) \in (0,T_M)\times\bar{\Omega}$}.
$$
%Observe that $x_0$ is not on $\Gamma$. Indeed, suppose that $x_0\in\Gamma$ then by maximum principle we have that
%$$
%e^{t\Delta_N}c_\varepsilon(t_0,x)> e^{t\Delta_N}c_\varepsilon(t_0,x_0),\quad (t,x)\in(0,\infty)\times\Omega
%$$
%Consequently from the definition of $\mathcal{K}_s$ we have that
%$
%\mathcal{K}_sc_\varepsilon(t_0,x)>\mathcal{K}_sc_\varepsilon(t_0,x_0),
%$
%consequently $\nabla\mathcal{K}_sc_\varepsilon(t_0,x_0)\cdot\nu>0$, which is a contradiction from the boundary condition in \eqref{eq_regulation}.
We claim that $c_\varepsilon(t_0,x_0)\geq0$, which implies that $c_\varepsilon$ is non-negative. 
Indeed, suppose that, $c_\varepsilon(t_0,x_0)<0$,
and evaluating $(t_0,x_0)$ in
the second equation \eqref{eq_regulation}, we obtain
\begin{equation}
(-\Delta_N)^{1-s}c_\varepsilon(t_0,x_0)+c_\varepsilon(t_0,x_0)=u_\varepsilon(t_0,x_0). 
\label{eq_value}
\end{equation}
Now, we recall that 
\begin{equation}
(-\Delta_N)^{1-s}c_\varepsilon(t_0,x_0)=\dfrac{1}{\Gamma(s-1)}\int^{\infty}_0
\big( e^{t\Delta_N}c_\varepsilon(t_0,x_0)-c_\varepsilon(t_0,x_0) \big) \dfrac{dt}{t^{2-s}},
\label{eq:DeltaK}
\end{equation}
where $\Gamma(s-1)< 0$ $(s< 1)$ and $v(t,x)= e^{t\Delta_N}c_\varepsilon(t_0,x)$ is the weak solution of the IBVP
$$
\left\lbrace\begin{aligned}
   \partial_t v-\Delta v&=0,\;\;\;\;\;\;\;\;\;\mbox{in}\;(0,\infty)\times\Omega,
   \\[3pt]
   \nabla v\cdot\nu&= 0,\;\;\;\;\;\;\;\;\;\mbox{on}\;[0,\infty)\times\Gamma,
   \\[3pt]
   v(0,x)&=c_\varepsilon(t_0,x),\mbox{in}\;\Omega.
\end{aligned}\right.
$$
Hence applying the (weak) maximum principle (see \cite{PW}, Theorem 7), we get that the minimum occur on  the parabolic boundary of $(0,\infty)\times\Omega$, which comprises 
$\left\lbrace0\right\rbrace\times\Omega$ and 
$[0,\infty)\times\Gamma$. Moreover, the minimum of $v$ could not occur on $[0,\infty)\times\Gamma$ since $\nabla v\cdot\nu=0$. Consequently, we obtain
$$
 e^{t\Delta}c_\varepsilon(t_0,x)\,\geq\,
c_\varepsilon(t_0,x_0),
$$  
for all $(t,x)\in(0,\infty)\times\Omega$. Therefore from \eqref{eq:DeltaK} we deduce that, 
$$
(-\Delta_N)^{1-s}c_\varepsilon(t_0,x_0) \leq 0,
$$
this together with \eqref{eq_value} and $c_\varepsilon(t_0,x_0)<0$ implies $u_\varepsilon(t_0,x_0)< 0$, which is a contradiction, hence $0\leq c_\varepsilon$. By an analogous argument, we obtain that $ c_\varepsilon(t,x)\leq 1$, which finish the proof.
\end{proof}
%%%

\begin{theorem}[Global Existence]\label{th_global} Given $u_0 \in L^{\infty}(\Omega)$ satisfying \eqref{ID},
then there exists a unique solution 
$$
    (u_\varepsilon, c_\varepsilon) \in W(\infty) \cap L^{\infty}((0,\infty) \times \Omega) 
    \times L^{\infty}\big((0,\infty); D\big((-\Delta_N)^{1-s}\big)\big)
$$    
of the problem \eqref{eq_regulation}, and it satisfies the following uniform bounds
\begin{eqnarray}
0 \leq u_{\varepsilon}(t, x) \leq 1,\quad  \text { a.e. in }[0,\infty) \times \Omega \label{Pro_eq_bound_u},
\\[5pt]
0 \leq c_{\varepsilon}(t, x) \leq 1,\quad  \text { a.e. in }[0,\infty) \times \Omega \label{Pro_eq_bound_c}. 
\end{eqnarray}
\end{theorem}
\begin{proof}
First, from Theorem \ref{th_local} there exists $T_M> 0$ and a unique solution $(u_\varepsilon, c_\varepsilon)$ 
of the problem \eqref{eq_regulation} posed in $(0,T_M)\times\Omega$. 
Since for almost all $x \in \Omega$, $0\leq u_0(x) \leq 1$, by Proposition \ref{pro_bound_u_c}, we have 
$0 \leq u_{\varepsilon}(t, x) \leq 1$ and $0 \leq c_{\varepsilon}(t, x) \leq 1$ for a.e. 
$(t,x) \in [0,T_M) \times \Omega$. 

\medskip
 \underline {Claim:} The maximal existence time $T_M=+\infty$.
 
 \medskip
Proof of Claim: Indeed, let us suppose that $T_M<+\infty$. Therefore, applying 
Remark \ref{rem_TM} there exists $t_j\to T_M^{-}$ as $j\to\infty$ such that
$$
\lim_{j\to\infty}\|u_\varepsilon(t_j,\cdot)\|_{L^{\infty}(\Omega)}= +\infty,
$$
which is a contradiction, thus $T_M= + \infty$. 
\end{proof}

%%%%%%%%%%%%%%%%%%%%%%%
\subsection{Perturbed problem estimates}
%%%%%%%%%%%%%%%%%%%%%%%

The aim of this section is to investigate some important properties of the global solution $(u_\varepsilon, c_\varepsilon)$ 
of the problem \eqref{eq_regulation} as given by Theorem \ref{th_global}. Henceforth, we consider that $(u_\varepsilon, c_\varepsilon)$ 
satisfies the partial differential equations of the problem \eqref{eq_regulation} in the strong sense, (see item 1 in the proof of
Proposition \ref{pro_bound_u_c}). 
 
\begin{lemma}
\label{Th_regu_entropi} Let $\left(u_{\varepsilon}, c_{\varepsilon}\right)$ be the unique solution of the problem \eqref{eq_regulation}. 
Then for any entropy $\eta \in \mathcal{C}^{2}$, 
\begin{equation}
\label{eq_entropy_perturbado}
\begin{aligned}
\frac{\partial}{\partial t} \eta\left(u_{\varepsilon}\right)+\operatorname{div}\left( q\left(u_{\varepsilon}\right)\nabla\mathcal{K}_s c_{\varepsilon}\right)&
+\left(u_{\varepsilon}- c_{\varepsilon}\right)\left[q-g \eta^{\prime}\right]\left(u_{\varepsilon}\right)
\\[5pt]
&+\varepsilon (-\Delta) \eta\left(u_{\varepsilon}\right) = -\varepsilon\eta^{\prime\prime}|\nabla u_{\varepsilon}|^2 \leq 0. 
\end{aligned}
\end{equation}
\end{lemma}

\begin{proof} First, multiplying the first equation in \eqref{eq_regulation}
by $\eta'(u_{\varepsilon})$ we obtain
\begin{equation}
\label{eq_entropy_perturbado_proof}
\partial_t \eta(u_{\varepsilon})+\varepsilon(-\Delta)\eta(u_{\varepsilon})
+\eta^{\prime}(u_{\varepsilon})\operatorname{div}\left(g(u_{\varepsilon})\nabla\mathcal{K}_s c_{\varepsilon} \right)
= -\varepsilon\eta^{\prime\prime}|\nabla u_{\varepsilon}|^2.
\end{equation}
On the other hand, one observes that 
$$
\begin{aligned}
\dive(q(u_{\varepsilon})\nabla\mathcal{K}_sc_{\varepsilon})&=& &q^{\prime}(u_{\varepsilon})\nabla u_{\varepsilon}\cdot\nabla\mathcal{K}_sc_{\varepsilon}-q(u_{\varepsilon})(-\Delta)^{1-s}c_{\varepsilon}
\\[5pt]
&=& &\eta'(u_{\varepsilon})g'(u_{\varepsilon})\nabla u_{\varepsilon}\cdot\nabla\mathcal{K}_sc_{\varepsilon}-q(u_{\varepsilon})(-\Delta)^{1-s}c_{\varepsilon},
\end{aligned}
$$
where we have used $q'(u)=\eta'(u)g'(u)$. Then from the second equation in \eqref{eq_regulation}, we have
$$
\begin{aligned}
\dive(q(u_{\varepsilon})\nabla\mathcal{K}_su_{\varepsilon})&=& &\eta'(u_{\varepsilon})\dive(g(u_{\varepsilon})\nabla\mathcal{K}_sc_{\varepsilon})+\left(\eta'(u_{\varepsilon})g(u_{\varepsilon})-q(u_{\varepsilon})\right)(-\Delta)^{1-s}c_{\varepsilon}
\\[5pt]
&=& &\eta'(u_{\varepsilon})\dive(g(u_{\varepsilon})\nabla\mathcal{K}_sc_{\varepsilon})-(u_{\varepsilon}-c_{\varepsilon})[q-g\eta'](u_{\varepsilon}),
\end{aligned}
$$
which substituting in \eqref{eq_entropy_perturbado_proof} concludes the proof.
\end{proof}
\begin{lemma}[Mass conservation]\label{Lem_mass_conser}
For almost all $t> 0$, 
$$
\int_{\Omega}c_{\varepsilon}(t,x) dx=  \int_{\Omega}u_{\varepsilon}(t,x)dx=\int_{\Omega}u_{0,\varepsilon}(x) dx.
$$
\end{lemma}
\begin{proof}
It is enough to integrate \eqref{eq_regulation} (first and second equations) over $\Omega$.
\end{proof}

\begin{proposition}
\label{Pro_grad_u} 
For any $T>0$, 
\begin{equation}
\label{EPSGRAD}
    \iint_{\Omega_T}|\sqrt{\varepsilon} \ \nabla u_{\varepsilon}(t,x)|^2\,dx\,dt\,\leq C,
\end{equation}
where $C= C(T, |\Omega|, \|u_0\|_{L^\infty})$ is a positive constant.
\end{proposition}
\begin{proof}
Let us multiply \eqref{eq_regulation} (first equation) by $u_{\varepsilon}$ and integrating by parts over $\Omega_T$, we obtain
\begin{equation}
\label{Proof_eq_1}
\begin{aligned}
\frac{1}{2}\int_{\Omega}u_{\varepsilon}(T)^2 dx
&+\varepsilon \iint_{\Omega_T} |\nabla u_{\varepsilon}(t,x)|^2 dxdt 
\\[5pt]
&=\frac{1}{2}\int_\Omega u_{0,\varepsilon}(x)^2 dx +\iint_{\Omega_T}g(u_{\varepsilon})\nabla\mathcal{K}_sc_{\varepsilon}\cdot\nabla u_{\varepsilon}(t,x) dxdt.
\end{aligned}
\end{equation}
Now, multiplying \eqref{eq_regulation} (second equation) by $\frac{u_\varepsilon^2}{2}-\frac{u_\varepsilon^3}{3}$, 
and integrating by parts we have 
$$
\begin{aligned}
\int_{\Omega}g(u_\varepsilon)\nabla\mathcal{K}_sc_\varepsilon\cdot\nabla u_\varepsilon(t,x) dx 
&=\,\int_{\Omega}  ((u_\varepsilon-u_\varepsilon^2)\nabla\mathcal{K}_sc_\varepsilon\cdot\nabla u_\varepsilon)(t,x) dx
\\[5pt]
&=\,\int_{\Omega} (\frac{u_\varepsilon^3}{2}-\frac{u_\varepsilon^4}{3}-\frac{c_\varepsilon\,u_\varepsilon^2}{2}+\frac{c_\varepsilon\,u_\varepsilon^3}{3})(t,x) dx.
\end{aligned}
$$
This, together with \eqref{Proof_eq_1} implies that 
$$
\begin{aligned}
\varepsilon\iint_{\Omega_T}|\nabla u_\varepsilon(t,x)|^2 dx dt &\leq \frac{1}{2}\int_\Omega u_{0,\varepsilon}(x)^2 dx 
\\[5pt]
&+\iint_{\Omega_T} (\frac{u_\varepsilon^3}{2}-\frac{u_\varepsilon^4}{3}-\frac{c_\varepsilon\,u_\varepsilon^2}{2}+\frac{c_\varepsilon\,u_\varepsilon^3}{3})(t,x) dx dt,
\end{aligned}
$$
which proves the assertion. 
\end{proof}

\begin{proposition}
\label{Lem_c_estima}
For each $T>0$ and $s \in (0,1)$, 
there exist positive constants $C_{1}= C_1(|\Omega|), C_2= C_{2}(T, |\Omega|)$, 
such that
\begin{eqnarray}
\left\| c_{\varepsilon}\right\|_{L^{\infty}\left((0, \infty); D((-\Delta_N)^{1-s})\right)} &\leq& C_{1}, \label{Pro_eq_bound_c_1-s}
\\[7pt]
\left\|\partial_{t} c_{\varepsilon}\right\|_{L^2((0,T);(-\Delta_N)^{s}(H^1(\Omega)))} &\leq& C_{2}.\label{Pro_eq_bound_c_t}
\end{eqnarray}
\end{proposition}

\begin{proof}
1. In order to show the first estimate, we recall from \eqref{Eq_c_epsilon} that   
$$
c_\varepsilon=\big(I_d+(-\Delta_N)^{1-s}\big)^{-1}u_\varepsilon 
= \sum_{k=0}^{\infty}\frac{1}{1+\lambda_k^{1-s}}\,\langle u_\varepsilon, \varphi_k\rangle\,\varphi_k,
$$
hence we obtain $c_\varepsilon \in L^{\infty}((0,\infty); D\big((-\Delta_N)^{1-s}\big))$. Moreover, for almost all $t>0$
$$
\|\,(-\Delta_N)^{1-s}c_\varepsilon(t)\|_{L^2(\Omega)}^2
= \sum_{k=0}^{\infty}\frac{\lambda_k^{2(1-s)}}{(1+\lambda_k^{1-s})^2}|\langle u_\varepsilon(t), \varphi_k\rangle|^2
\leq \sum_{k=0}^{\infty}|\langle u_\varepsilon(t), \varphi_k\rangle|^2<\infty.
$$

\medskip
2. To show \eqref{Pro_eq_bound_c_t}, we first observe that from \eqref{Eq_c_epsilon}, we have for almost all $t> 0$
$$
\begin{aligned}
\partial_t c_\varepsilon(t)&=\big(I+(-\Delta_N)^{1-s}\big)^{-1}\partial_tu_\varepsilon(t)
\\[5pt]
&=\big(I+(-\Delta_N)^{1-s}\big)^{-1}\big(\,-\dive(g(u_\varepsilon(t))\nabla\mathcal{K}_sc_\varepsilon(t))+\varepsilon \Delta u_\varepsilon(t) \big)
\end{aligned}
$$
in distributional sense. 

\medskip
On the other hand, from Proposition \ref{Pro_grad_u} together with \eqref{Pro_eq_bound_c_1-s}, we have that
$$
-\dive(g(u_\varepsilon)\nabla\mathcal{K}_sc_\varepsilon)+\varepsilon \Delta u_\varepsilon
$$
is uniformly bounded (with respect to $\varepsilon> 0)$ in $L^2((0,T);(H^1(\Omega))^{\star})$. 
One recalls that $(H^1(\Omega))^{\star}= (-\Delta_N)(H^1(\Omega))$. 
Finally, applying Proposition 2.1 we get the result. 
\end{proof}
%%%%%%%%%%%%%%%%%%%%%%%%%%%%%%%%%%%%%%%%%%%%%%%%%
%%%%%%%%%%%%%%%%%%%%%%%%%%%%%%%%%%%%%%%%%%%%%%%%%
%%%%%%%%%%%%%%%%%%%%%%%%%%%%%
%\begin{remark}
%Observe o seguinte
%$$
%\begin{aligned}
%H^1(\Omega)\,\subset\,& (-\Delta_N)^s(H^1(\Omega))\,\subset\, L^2(\Omega),\quad & &\text{ for } 0<s<1/2,
%\\[7pt]
%&(-\Delta_N)^s(H^1(\Omega))=L^2(\Omega),\quad& & s=1/2
%\\[7pt]
%L^2(\Omega)\,\subset\,& (-\Delta_N)^s(H^1(\Omega))\,\subset\, (H^1(\Omega))^{\star},\quad & &\text{ for } 1/2<s<1
%\end{aligned}
%$$
%on the other hand, we also have
%$$
%\begin{aligned}
% (-\Delta_N)^s(H^1(\Omega))&=D\big((-\Delta)^{1/2-s}\big),\quad & &\text{ for } 0<s<1/2,
% \\[7pt]
%(-\Delta_N)^s(H^1(\Omega))&=L^2(\Omega),\quad& & s=1/2
%\\[7pt]
% (-\Delta_N)^s(H^1(\Omega))&=\big(D\big((-\Delta)^{s-1/2}\big)\big)^{\star},\quad & &\text{ for } 1/2<s<1
%\end{aligned}
%$$
%\end{remark}
%%%%%%%%%%%%%%%%%%%%%%%%%%%%%%%%%%%%%%%%%%%%%%%%%
\begin{comment}
\begin{corollary}
For each $\varepsilon> 0$, $c_\varepsilon\in L^{\infty}((0,\infty);H^{1-s}(\Omega))$ and satisfies
$$
\left\| c_{\varepsilon}\right\|_{L^{\infty}\left((0, \infty) ; H^{1-s}(\Omega))\right)} \leq C_{1}, 
%\label{Pro_eq_bound_c_1-s}
$$
where $C_1$ is a positive constant.
\end{corollary}
\begin{proof}
First, from item $(4)$ in Proposition \ref{THMCARAC} and Proposition \ref{Lem_c_estima}, it follows that 
$$
\left\| c_{\varepsilon}\right\|_{L^{\infty}\left(0, \infty ; D((-\Delta_N)^{(1-s)/2})\right)} \leq C_{1},
$$
this together with Lemma \ref{Lem_mass_conser} and Proposition \ref{represendoDominio}  complete the proof.
\end{proof}
\end{comment} 

%%%%%%%%%%%%%%%%%%%%
\section{Proof of Main Theorem}
\label{LP}
%%%%%%%%%%%%%%%%%%%%

In this section we prove Theorem \ref{ETGS}, and to this end 
we are mostly concerned to pass to the limit in \eqref{eq_entropy_perturbado} as $\varepsilon \rightarrow 0$. 
More precisely, we write \eqref{eq_entropy_perturbado} in weak sense (using the entropy pair $\mathbf{F}(u, v)$), and jointly with 
\eqref{eq_regulation} (second equation) we obtain \eqref{CLE}, \eqref{EE} respectively, after pass to the limit as $\varepsilon \to 0$. 
Since  \eqref{eq_entropy_perturbado} has non-linear terms,
the uniform estimates on the sequence $\{u_\varepsilon\}$ are
not sufficient to take the limit transition on $\varepsilon $ as it goes to $%
0$. In fact, we need strong convergence, 
and as usual for scalar conservation laws we apply 
(following closer Chemetov, Neves \cite{Chemetov-Neves} and also Perthame, Dalibard \cite{perth-dali-tams}) the Kinetic Theory.

\medskip
First, we take in Lemma \ref{Th_regu_entropi}, 
$\eta(u)=|u-k|^{+}$, $q(u)=\sgn^{+}(u-k)(g(u)-g(k))$ with $k \in \mathbb{R}$. Then, we obtain from equation \eqref {eq_entropy_perturbado} and system \eqref{eq_regulation}, 
for each test function $\phi\in C^{\infty}_0\left((-\infty,T)\times\R^n\right)$, (for simplicity of exposition and abuse of notation we may have $T= \infty$), 
\begin{equation}
\label{eq_int_entro_+}
\begin{aligned}
\iint_{\Omega_T} &(\eta(u_\varepsilon)\phi_t\,+\,q(u_\varepsilon)\nabla\mathcal{K}_sc_\varepsilon\cdot\nabla\phi-\varepsilon\nabla\eta(u_\varepsilon)\cdot\nabla\phi) \ dx dt
\\[5pt]
&+ \int_{\Omega}\eta(u_{0,\varepsilon})\phi(0) dx + \iint_{\Omega_T}\left(u_{\varepsilon}-c_{\varepsilon}\right)\eta^{\prime}(u_{\varepsilon})g(k)\phi \ dx dt 
\\[5pt]
& =m_{\varepsilon}^+(\phi),
\end{aligned}
\end{equation}
where we have used $\left[q-g \eta^{\prime}\right]\left(u_{\varepsilon}\right)=-\eta^{\prime}(u_{\varepsilon})g(k)$,
and $m_{\varepsilon}^+$ is a real non-negative Radon measure, defined by
\begin{equation}
\label{MedidaPosit}
m_{\varepsilon}^+(\phi):= \iint_{\Omega_T} \varepsilon  \ \eta^{\prime\prime}(u_\varepsilon) \ |\nabla u_{\varepsilon}|^2 \phi \ dx dt.
\end{equation}

\medskip
Now, we differentiate in the distributional sense equation \eqref {eq_entropy_perturbado} with respect to $k$, 
hence we obtain (as now a standard procedure in the kinetic theory) the following transport like equation 
\begin{equation}
\frac{\partial f_{\varepsilon}}{\partial t}+\left(k-c_{\varepsilon}\right) g(k) \frac{\partial f_{\varepsilon}}{\partial k}
+g^{\prime}(k) \nabla  \mathcal{K}_s c_{\varepsilon} \cdot \nabla f_{\varepsilon}+\varepsilon (-\Delta) f_{\varepsilon}=\partial_k m_{\varepsilon}^+,
\label{eq_kinetic_+}
\end{equation}
where  $f_\varepsilon(t,x,k):=\sgn^{+}(u_{\varepsilon}(t,x)-k)$.
Rigorously, from \eqref{eq_int_entro_+} we get that the function $f_\varepsilon(t,x,k)$ satisfies the following equation
\begin{equation}
\begin{aligned}
\iint_{\Omega_T} & \big\lbrace\int_k^1f_\varepsilon(t,x,v)[\phi_t+g'(v)\nabla\mathcal{K}_sc_\varepsilon\cdot\nabla\phi] dv-\varepsilon\nabla\eta(u_\varepsilon)\cdot\nabla\phi \big\rbrace dx dt
\\[5pt]
&+ \int_{\Omega}|u_{0,\varepsilon}-k|^{+}\phi(0) \ dx 
\\[5pt]
&+ \iint_{\Omega_T}\left(u_{\varepsilon}-c_{\varepsilon}\right)g(k)f_\varepsilon\,\phi \ dx dt =m_{\varepsilon}^+(\phi)\geq 0,
\end{aligned}\label{eq_int_kinetic_+}
\end{equation}
for all nonnegative function $\phi\in C^{\infty}_0\left((-\infty,T)\times\R^n\right)$.
Furthermore, it follows that for any function $G\in C^{1}([0,1])$, with $G(0)= 0$,
\begin{eqnarray}
G(u_{\varepsilon }) &=&\int_{0}^{1}G^{\prime }(v)f_{\varepsilon }(\cdot ,
\mathbf{\cdot },v)\; dv \quad \text{a.e. in $\Omega _{T},$}  \notag 
\\
0 &\leq &f_{\varepsilon} \leq 1 \quad \text{on $\Omega _{T}\times
\mathbb{R}$},\quad \text{$f_{\varepsilon }(t,x,k)
=\left\{ 
\begin{aligned} 1, &\quad \text{for $k \leq 0$}, 
\\ 
0, &\quad \text{for $k \geq 1$}, 
\end{aligned}\right. $}  \notag 
\\
\partial _{k} f_{\varepsilon } &\leq &0\qquad \text{in ${\mathcal{D}}%
^{\prime }(\Omega _{T}\times \mathbb{R})$}. 
 \label{u-f}
\end{eqnarray}
From \eqref{Pro_eq_bound_u}, \eqref{EPSGRAD}, \eqref{MedidaPosit} and the 
Riesz Representation Theorem, we get that 
$m_{\varepsilon}^+$ is a real positive Radon measure, defined on $\overline{%
\Omega }_{T}\times \mathbb{R},$ and
\begin{align}
& m_{\varepsilon}^+(\cdot ,\cdot ,k)=0 \quad \text{for any } k> 1\quad
\text{on }\overline{\Omega}_{T},\;\text{and for any } |\phi |\leq 1 \notag 
\\[5pt]
& \iint_{\overline{\Omega }_{T}\times \mathbb{R}}\phi \;dm_{\varepsilon}^{+}\leq C \quad \text{continuously
\ on }\overline{\Omega }_{T}\times \mathbb{R},  \label{m-property}
\end{align}%
where $C$ is a positive constant independent of $\varepsilon .$

\medskip
Similarly to derive the inequality \eqref{eq_int_kinetic_+}, 
we now consider $\eta(u)=|u-k|^{-}$, and $q(u)=\sgn^{-}(u-k)(g(u)-g(k))$. Then, we obtain
\begin{equation}
\begin{aligned}
\frac{\partial }{\partial t}(1-f_{\varepsilon})+\left(k-c_{\varepsilon}\right) g(k) \frac{\partial}{\partial k}(1-&f_{\varepsilon})+g^{\prime}(k) \nabla  \mathcal{K}_s c_{\varepsilon} \cdot \nabla (1-f_{\varepsilon})
\\[5pt]
&+\varepsilon (-\Delta)(1-f_{\varepsilon})=-\partial_k m_{\varepsilon}^-
\end{aligned}\label{eq_kinetic_-}
\end{equation}
in the distribution sense, and for all nonnegative function $\phi\in C^{\infty}_0\left((-\infty,T)\times\R^n\right)$, we have the following identity
\begin{equation}
\begin{aligned}
\iint_{\Omega_T} &\big\lbrace\int_0^k(1-f_\varepsilon(t,x,v))[\phi_t+g'(v)\nabla\mathcal{K}_sc_\varepsilon\cdot\nabla\phi] dv-\varepsilon\nabla\eta(u_\varepsilon)\cdot\nabla\phi \big\rbrace dx dt
\\[5pt]
&+\int_{\Omega}|u_{0,\varepsilon}-k|^{-}\phi(0) \ dx 
\\[5pt]
&+\iint_{\Omega_T}\left(u_{\varepsilon}-c_{\varepsilon}\right)g(k)(1-f_\varepsilon)\,\phi \ dx dt= m_{\varepsilon}^-(\phi)\geq 0,
\end{aligned}\label{eq_int_kinetic_-}
\end{equation}
where $m_{\varepsilon}^-$ is defined in the same way by \eqref{MedidaPosit}. 
Moreover, the real positive Radon measure $m_{\varepsilon}^{- },$
defined on $\overline{\Omega }_{T}\times \mathbb{R}$,\ satisfies the
following properties
\begin{align}
m_{\varepsilon}^{- }(\cdot ,\cdot ,k)& =0\quad \text{for any }k<0\text{
on }\overline{\Omega }_{T}, \;\text{and for any } |\phi |\leq 1  \notag 
\\[5pt]
\iint_{\overline{\Omega }_{T}\times \mathbb{R}}\phi \;dm_{\varepsilon}^{- }&\leq C \text{\quad continuously
on }\overline{\Omega }_{T}\times \mathbb{R}.  \label{m-property2}
\end{align}

\medskip
At this point, let us study the convergence of the sequence $\{c_\varepsilon\}$. 
First, from the uniform estimate \eqref{Pro_eq_bound_c}, there exists a function 
$c \in L^\infty((0,\infty) \times \Omega)$, such that (passing to a subsequence)
\begin{equation}
\label{C1}
c_\varepsilon \rightharpoonup c \quad \text{weakly$-\star$ in }L^{\infty}_\loc((0,\infty) ; L^{\infty}(\Omega)). 
\end{equation}
\begin{comment} 
Moreover, since $L^{\infty}_\loc((0,\infty); L^{\infty}(\Omega)) \subset L^{\infty}_\loc((0,\infty); L^{2}(\Omega))$ 
it follows from Lemma \ref{Lem_mass_conser} that
\begin{equation}
\label{C2}
   %\int_\Omega c_\varepsilon(t) \ dx \to 
   \int_\Omega c(t) \ dx= 0???????(NAO). 
\end{equation}
\end{comment}
%
Although, we need more than convergence in averages for the sequence $\{c_\varepsilon\}$. 
Then, we have the following 
\begin{proposition} 
\label{PROPREGC}
Given $s \in (0,1)$, there exists $c \in L^{\infty}\big((0,\infty);D((-\Delta_N)^{1-s})\big)$, such that 
\begin{equation}
\label{C3}
     c_\varepsilon \to c \quad {\rm strongly \ in \ }  L^2_{\rm loc}((0,\infty); D\big((-\Delta_N)^{(1-s)/2}\big).
\end{equation}
\end{proposition}
\begin{proof} 
First, from the uniform estimate \eqref{Pro_eq_bound_c_1-s}, and passing to a convenient subsequence, there exists a 
function $c \in L^{\infty}\big((0,\infty);D((-\Delta_N)^{1-s})\big)$,
such that
$$
   c_\varepsilon \rightharpoonup c \quad \text{ weakly in } L^{2}_\loc \big((0,\infty);D((-\Delta_N)^{1-s})\big).
$$

\medskip
Now, due to \eqref{Pro_eq_bound_c_t} $\partial_tc_\varepsilon$ is uniform bounded in $L^2_\loc((0,\infty);(-\Delta_N)^s(H^1(\Omega)))$. 
Finally, from definitions \eqref{2equalrepredominio}, \eqref{equalrepredominio} and applying the Aubin-Lions' Theorem, we get \eqref{C3}. 
\end{proof}

\medskip
Therefore, in view of the above results (passing to subsequences still denoted by $\varepsilon$), 
there exist functions $u \in L^{\infty}((0, \infty) \times \Omega)$, $f \in L^{\infty}((0, \infty) \times \Omega \times \mathbb{R})$, 
$c \in L^{\infty}\left((0,\infty) ; D((-\Delta_N)^{1-s})\right)$, and non-negative measures $m^{\pm}= m^{\pm}(t, x, k)$, such that (locally in time),
\begin{eqnarray}
u_{\varepsilon} \rightharpoonup u, &\quad & \text{weakly$-\star$ in } L^{\infty},
\nonumber
\\[5pt]
f_{\varepsilon} \rightharpoonup f, &\quad & \text{weakly$-\star$ in } L^{\infty},
\nonumber
\\[5pt]
m_{\varepsilon}^\pm \rightharpoonup m^{\pm}, &\quad & \text{weakly in } \mathcal{M}^+, 
\nonumber
\\[5pt]
c_{\varepsilon} \rightarrow c, &\quad & \text {strongly in } L^{2}_\loc\left((0,\infty); L^2(\Omega)\right),
\nonumber
\\[5pt]
\nabla\mathcal{K}_sc_{\varepsilon} \rightarrow \nabla\mathcal{K}_sc, &\quad & \text {strongly in } L^{2}_\loc\left((0,\infty); \mathbf{L}^2(\Omega)\right).
\nonumber
\end{eqnarray}
The above convergences are enough to pass to the limit (as $\varepsilon> 0$ goes to zero) 
in the second equation of the system \eqref{eq_regulation}, that is to say 
\begin{equation}
\label{C4}
(-\Delta_N)^{1-s} c + c = u
\end{equation}
also in equations \eqref{eq_kinetic_+} and \eqref{eq_kinetic_-}, the only exception is the term 
$g^{\prime}(k) \nabla\mathcal{K}_s c_{\varepsilon} \cdot \nabla f_{\varepsilon}$,
which yields an extra effort. First, we can write
$$
\begin{aligned}
g^{\prime}(k) \nabla\mathcal{K}_s c_{\varepsilon} \cdot \nabla f_{\varepsilon} &=\dive\left(g^{\prime}(k)f_{\varepsilon} \nabla\mathcal{K}_s c_{\varepsilon} \right)+g^{\prime}(k) (-\Delta_N)^{1-s} c_{\varepsilon} f_{\varepsilon}
 \\[5pt]
&=\dive\left(g^{\prime}(k)f_{\varepsilon} \nabla\mathcal{K}_s c_{\varepsilon} \right)+\left(u_{\varepsilon}-c_{\varepsilon}\right) g^{\prime}(k) f_{\varepsilon}.
\end{aligned}
$$
Moreover, we have in the sense of distributions as $\varepsilon \to 0$, 
\begin{eqnarray}
\dive\left(g^{\prime}(k)f_{\varepsilon} \nabla\mathcal{K}_sc_{\varepsilon} \right) &\rightarrow& \dive\left(g^{\prime}(k)f \nabla\mathcal{K}_sc \right),
\nonumber
\\[5pt]
c_{\varepsilon} \ g^{\prime}(k) f_{\varepsilon} &\rightarrow& c \ g^{\prime}(k) f. 
\nonumber
%\\[5pt]
%\langle u_{\varepsilon}\rangle g^{\prime}(k) f_{\varepsilon}=\langle u_{0,\varepsilon}\rangle g^{\prime}(k) f_{\varepsilon} &\rightarrow& \langle u_0\rangle g^{\prime}(k) f. \nonumber
\end{eqnarray}
Although, from the moment,
we cannot assert that the weak limit of $u_{\varepsilon} f_{\varepsilon}$ is $u f$. 
However, we know that $\left\{u_{\varepsilon} f_{\varepsilon}\right\}_{\varepsilon>0}$ is 
uniformly bounded in $L^{\infty}((0,\infty) \times \Omega \times\R)$.
Then, extracting a further subsequence (if necessary), there exists a function $\rho= \rho(t, x,k)
\in L^{\infty}((0,\infty) \times \Omega \times\R)$, such that
(locally in time) 
\begin{equation}
u_{\varepsilon} f_{\varepsilon} \rightharpoonup \rho, \quad \text{weakly$-\star$ in } L^{\infty}.
\label{Eq_rho}
\end{equation}
%\begin{comment}
\begin{remark}
\label{rem_rho}
Thanks to the definition of $f_\varepsilon$, together with \eqref{Eq_rho}, one observes that 
$$
\rho(t,x,k)=
\begin{cases}
0,\quad \quad\quad\text{ when } k\geq  1,
\\[5pt]
u(t,x),\quad \text{ when } k\leq 0,
\end{cases}
$$
almost everywhere in $(0,\infty) \times \Omega \times \R$.
\end{remark}
%\end{comment}
Consequently, we obtain in distribution sense 
$$
g^{\prime}(k) \nabla\mathcal{K}_sc^{\varepsilon} \cdot \nabla f^{\varepsilon} \rightharpoonup g^{\prime}(k) \nabla\mathcal{K}_sc \cdot \nabla f+g^{\prime}(k)(\rho-u f),
$$
where we have used \eqref{C4}.
Therefore, from \eqref{eq_int_kinetic_+} and  \eqref{eq_int_kinetic_-}, it follows respectively that, 
 for any nonnegative function $\phi\in C^{\infty}_0\left((-\infty,T)\times\R^{n+1}\right)$,
\begin{equation}
\begin{aligned}
\iint_{\Omega_T} & \int_k^1f(t,x,v)[\phi_t+g'(v)\nabla\mathcal{K}_sc\cdot\nabla\phi] dv \ dx dt
\\[5pt]
& +\iint_{\Omega_T}(\rho-uf) g(k)\phi \ dx dt
+\iint_{\Omega_T}\left(u-c\right)g(k)f \phi \ dx dt
\\[5pt]
&+\int_{\Omega}|u_{0}-k|^{+}\phi(0)  \ dx= m^+(\phi) =: \iint_{\Omega_T} m^+(t,x,k) \phi \ dxdt  \geq 0,
\end{aligned}\label{eq_int_f_kinetic_+}
\end{equation}

and
\begin{equation}
\begin{aligned}
\iint_{\Omega_T} &\int_0^k(1-f(t,x,v))[\phi_t+g'(v)\nabla\mathcal{K}_sc\cdot\nabla\phi] dv \ dx dt
\\[5pt]
& -\iint_{\Omega_T}(\rho-uf) g(k)\phi \ dx dt  
+\iint_{\Omega_T}\left(u-c\right)g(k)(1-f)\,\phi \ dx dt
\\[5pt]
& +\int_{\Omega}|u_{0}-k|^{-}\phi(0) \ dx = m^-(\phi) =: \iint_{\Omega_T} m^-(t,x,k) \phi \ dxdt  \geq 0.
\end{aligned}\label{eq_int_f_kinetic_-}
\end{equation}
Moreover, for any function $G\in C^{1}([0,1])$, with $G(0)= 0$, it follows that
\begin{eqnarray}
G(u) &=&\int_{0}^{1}G^{\prime }(v)f(\cdot ,\mathbf{\cdot },v)\;dv
\qquad \text{ a.e. in $\Omega _{T},$} 
 \notag 
\\
0 &\leqslant &f\leqslant 1\quad \text{on $\Omega _{T}\times \mathbb{R}$}%
,\quad \text{$f(t,x,k)=\left\{ \begin{aligned} 1, &\quad \text{for
$k \leq 0$}, \\ 0, &\quad \text{for $k \geq 1$}, \end{aligned}\right. $}
\notag 
\\
\partial _{k}f &\leq &0 \quad \text{in distribution sense,}  \label{for u}
\end{eqnarray}%
also we have 
\begin{align}
\iint_{\overline{\Omega }_{T}\times \mathbb{R}}\phi \;dm^{\pm }& \leq
C \quad \text{ for any }|\phi |\leq 1,  
\notag 
\\
m^{+}(\cdot ,\cdot ,k)& = 0 \quad \text{for any }k>1\text{ and on }\overline{%
\Omega }_{T},  
\notag 
\\
m^{-}(\cdot ,\cdot ,k)& = 0 \quad \text{for any }k<0\text{ and on }\overline{%
\Omega }_{T},  
\notag 
\\
m^{\pm }(\cdot ,\cdot ,k)& \in C(\mathbb{R};\text{${\mathcal{M}}^{+}(%
\overline{\Omega }_{T}\times \mathbb{R})$}),  \label{for m}
\end{align}%
where the continuity of $m^{\pm }$ on $k$ \ follows from 
\eqref{eq_int_f_kinetic_+}, \eqref{eq_int_f_kinetic_-}.

\medskip 
Finally, we take $\phi =\partial_{k}\psi$ in \eqref{eq_int_f_kinetic_+} and \eqref{eq_int_f_kinetic_-}, with $\psi $ being a nonnegative function in $C_{0}^{\infty
}(\Omega _{T}\times \mathbb{R})$. Then integrating by parts on $k$, we obtain
that $f$ satisfies, respectively, the following transport equations (in distribution sense)
\begin{equation}
\frac{\partial}{\partial t} f+\mathbf{b}\cdot \nabla_{\!\!(k,x)} f
\,+\,g^{\prime}(k)(\rho-u f)\,=\,\partial_{k} m^{+},
\label{eq_+}
\end{equation}
and 
\begin{equation}
\frac{\partial}{\partial t}(1-f)+\mathbf{b} \cdot \nabla_{\!\!(k,x)} (1-f)
\,-\,g^{\prime}(k)(\rho-uf)=-\partial_{k} m^{-}.
\label{eq_-}
\end{equation}
Here $\nabla_{\!\!(k,x)} =\left(\frac{\partial}{\partial k}\, ,\,\nabla_x\right)$ and the vector field $\mathbf{b}:(0,T)\times\Omega \times [0,1] \to \R\times\R^n$,
called drift, is given by 
\begin{equation}
\label{DRIFTB}
\mathbf{b}(t,x,k)=\big( (k -c(t,x))g(k)\,,\,g^{\prime}(k)\nabla\mathcal{K}_sc(t,x) \big).
\end{equation}
Moreover, we have for $0< s \leq 1/2$
\begin{equation}
\label{renormalization} 
\begin{aligned}
&\mathbf{b}\in L^\infty\left((0,\infty);\mathbf{H}^1(\Omega \times [0,1] )\right), \quad \text{and}
\\[5pt]
 & \mathbf{b} \cdot \nabla_{\!\!(k,x)} f=  \dive_{\!(k,x)} \big(\mathbf{b} f \big) \quad \text{in distribution sense.}  
\end{aligned}    
\end{equation}

\begin{lemma}\label{Lem_f(1-f)} Let $\mathbf{b}$ be the drift vector field defined by \eqref{DRIFTB},
and $0< s \leq 1/2$. 
Then, the function $F= f (1-f)$ satisfies in the sense of distributions
\begin{equation}
\frac{\partial }{\partial t} F + \dive_{\!(k,x)} \big(\mathbf{b} F \big) + R \ (1-2f)\,\leq\,0,
\label{eq_f(1-f)}
\end{equation}
where $R:= g^{\prime}(k)(\rho-u f)$.
\end{lemma}
\begin{proof}
Under the conditions of the vector field $\mathbf{b}$ given by \eqref{renormalization}, 
we can apply the renormalization procedure which means that, the equations 
\eqref{eq_+} and \eqref{eq_-} are regularized on a parameter $\theta$, and respectively 
multiplied by $(1-f^\theta)$ and $f^\theta$, ($f^\theta$ being the regularization of $f$).   
Then, the obtained equations are added and, taking the limit as $\theta \to 0$ we obtain 
\eqref{eq_f(1-f)}, where the inequality follows from the following relation 
$$
   \int_{\R}\left((1-f^\theta) \partial_k (m^{+})^\theta\,-\,f^\theta \partial_k(m^-)^\theta\right)\,dk\,
   =\,\int_{\R}(m^{+}\,+\,m^{-})^\theta \partial_k f^\theta\,dk\,\leq\, 0,
$$ 
in view of \eqref{for u}, \eqref{for m}. Actually, we omit the details as now it is a standard procedure
 in the renormalization theory for transport equations. 
\end{proof}

Now, let us study the trace concept of $f$ at time $t=0$.  
\begin{proposition}
\label{pro_intial}
The function $f(t,x,k)$ has the trace $f_0(x,k)$ at time $t=0$, such that
$$
    f_0(x,k):=\ess\lim_{t\to0}f(t,x,k) \quad \text{almost everywhere in $\Omega\times\R$.}
$$
Moreover, $f_0= (f_0)^2$. 
\end{proposition}
\begin{proof}
1. First, let $\clg{E}$ be a countable dense subset of $C^{1}_0(\Omega)$. Then, 
for each $\zeta \in \clg{E}$ and $k\in\mathbb{Q}\cap[0,1]$, we define the following set of full measure in $(0,T)$,
$$
E_{\zeta,k}:= \Big\{ t \in(0,T) / \, t  \text{ is a Lebesgue point of }  I(t)= \int_{\Omega}\int_k^1f(t,x,v)\zeta(x)dvdx \},
$$
 and consider
$$
E:=\bigcap_{(\zeta,k)} E_{\zeta,k}, 
$$
where the intersection is taken over
$\clg{E}\times\left(\mathbb{Q}\cap[0,1]\right)$. Also $E$ is a set of full measure in $(0,T)$.

\medskip
2. To show the existence of the essential limit of $f(t,x,k)$, as $t$ goes to zero, we use the inequalities \eqref{eq_int_f_kinetic_+},
\eqref{eq_int_f_kinetic_-}. Indeed, we consider the test function $\phi(t,x,k)= \zeta_j(t) \psi(x)$, 
$\zeta_j(t)=H_j(t+t_0)-H_j(t-t_0)$ for any $t_0 \in E$ (fixed), and $j\geq 1$,
where $H_j$ is a standard regularization of the Heaviside function, and 
$\psi$ is a non-negative function which belongs to $\clg{E}$. 
Then, we have from \eqref{eq_int_f_kinetic_+} 
$$
\begin{aligned}
\int_0^{T} \int_{\Omega}\int_k^1 f(t,x,v) \,\zeta^{\prime}_j(t) \psi(x) \ dv dx dt 
&+ \int_{0}^T \Phi_k(t) \,\zeta_j(t) \ dt 
\\[5pt]
&+\int_{\Omega}|u_{0}-k|^{+} \ \psi(x) \ dx \geq 0,
\end{aligned}
$$
where
$$
\begin{aligned}
\Phi_k(t)= \int_{\Omega} &\Big( \int_k^1 f(t,x,v) \ g'(v) \nabla\mathcal{K}_sc \cdot \nabla \psi \ dv  
\\[5pt]
&+ \big( (\rho-uf) g(k)+ (u-c)\big) \ g(k) \ f(t,x,k) \ \psi(x) \Big) \ dx. 
\end{aligned}
$$
Passing to the limit in the above equation as $j \to \infty$, and taking into account that $t_0$ is Lebesque point of $I(t)$, we obtain
\begin{equation}
\label{INT100}
   \int_{\Omega} \psi(x) \left\lbrace-\int_{k}^1 f(t_0,x,v) \ dv + |u_0(x)-k|^{+}\right\rbrace dx 
   + \int^{t_0}_0 \Phi_k(t) \ dt \geq 0,
\end{equation} 
where we have used the Dominated Convergence Theorem. Since $t_0 \in E$ is arbitrary, and in view of the density of $\clg{E}$ 
in $L^1(\Omega)$, it follows from \eqref{INT100} that 
$$
    \ess\lim_{t\to 0}\int_{\Omega} \left\lbrace-\int_{k}^1f(t,x,v)dv  + |u_0(x)-k|^{+} \right\rbrace \psi(x) \ dx \geq 0
$$
for all non-negative $\psi \in L^1(\Omega)$, which implies for almost everywhere $x$ in $\Omega$,
$$
\ess\lim_{t\to 0}f(t,x,k)=0,\quad \text{ if } k> u_0(x). 
$$
Similarly, we obtain from \eqref{eq_int_f_kinetic_-}
$$
    \ess\lim_{t\to 0}\int_{\Omega}\left\lbrace- \int_{0}^k(1-f(t,x,v))dv  + |u_0(x)-k|^{-} \right\rbrace\ \psi(x) \ dx \geq 0
$$ 
for all non-negative $\psi \in L^1(\Omega)$, which implies for almost everywhere $x$ in $\Omega$,
$$
\ess\lim_{t\to 0}f(t,x,k)=1,\quad \text{ if } k<u_0(x).
$$ 
Therefore, the $\ess\lim_{t\to 0} f(t,x,k)$ exists, and in particular we have 
$$f_0(x,k)= \sgn^{+}(u_0(x)-k)$$
almost everywhere in $\Omega\times\R$,
which concludes the proof.
\end{proof}

One remarks that, since $f \in L^{\infty}((0,\infty) \times \Omega \times\R)$, it follows that
\begin{equation}
\label{TRACE0}
   \ess\lim_{t\to0^+} f(t,x,k)= \lim_{\delta\to0^+}\frac{1}{\delta} \int_0^{\delta} f(\tau,x,k) \,d\tau
\end{equation}
almost everywhere in $\Omega\times\R$. 

\medskip
Before we gain the strong convergence of $u_\epsilon$, which is obtained showing that,
$f= f^2$, which is to say, $F(t,x,k)= 0$ almost everywhere in $(0,\infty) \times \Omega \times \R$,
it remains to study the remainder term $R$, that is 
$$
R(t, x,k)= g^{\prime}(k)(\rho-u f)(t, x,k).
$$
In fact, this study has been done in \cite{perth-dali-tams}, and we recall here the main details 
with minor modifications. 
%for completeness of the  paper. 
First, from Remark \ref{rem_rho} we only need to obtain the formula for $\rho$, once $k \in (0,1)$, 
since $R \equiv 0$ for $k \leq 0$ and $k \geq 1$. 
Then, considering the test functions 
$\varphi_{1} \in \mathcal{C}_{0}^{\infty}((-\infty, T) \times \Omega)$, $\varphi_{2} \in \mathcal{C}_{0}^{\infty}(0,1)$, for any $T> 0$, 
we have
$$
\begin{aligned}
\iint_{\Omega_T} \int_0^1 & \rho(t, x, k) \varphi_{1}(t, x) \varphi_{2}^{\prime}(k) \,dk \ dx \ dt
\\[5pt]
=& \lim _{\varepsilon \rightarrow 0} \iint_{\Omega_T} \int_0^1 u_{\varepsilon}(t, x) f_{\varepsilon}(t, x, k)  \varphi_{1}(t, x)\varphi_{2}^{\prime}(k)\,dk \ dx\, dt  
\\[5pt]
=& \lim _{\varepsilon \rightarrow 0} \iint_{\Omega_T} u_{\varepsilon}(t, x) \varphi_{1}(t, x)\left\lbrace \int_0^1 \varphi_{2}^{\prime}(k)  f_{\varepsilon}(t, x, k) \,dk\right\rbrace \,dx \, dt
%\\[5pt]
%=& \lim _{\varepsilon \rightarrow 0} \int_{0}^{\infty} \int_{\Omega} u_{\varepsilon}(t, x) \varphi_{1}(t, x)\left\lbrace\int_0^1 \delta_{k=u_{\varepsilon}(t,x)} \varphi_{2}(k)\,dk\right\rbrace \,dt\, dx  
\\[5pt]
=& \lim _{\varepsilon \rightarrow 0} \iint_{\Omega_T} u_{\varepsilon}(t, x)\,\varphi_{2}\left(u_{\varepsilon}(t, x)\right) \varphi_{1}(t, x)\,dx \, dt
%\\[5pt]
%=& \lim _{\varepsilon \rightarrow 0} \int_{0}^{\infty} \int_{\Omega}  \varphi_{1}(t, x)\left\lbrace\int_\R \delta_{k=u_{\varepsilon}(t,x)}\,k\,\varphi_{2}(k)\,dk\right\rbrace \,dt\, dx  
\\[5pt]
=& \lim _{\varepsilon \rightarrow 0} \iint_{\Omega_T} \left\lbrace \int_0^1 \frac{d}{d k}\left(k\, \varphi_{2}(k)\right) f_{\varepsilon}(t, x, k) \ dk \right\rbrace \varphi_{1}(t, x)\, dx \, dt
\\[5pt]
=& \iint_{\Omega_T} \left\lbrace \int_0^1 \frac{d}{d k}\left(k\, \varphi_{2}(k)\right) f(t, x, k) \ dk \right\rbrace \varphi_{1}(t, x)\, dx \, dt,
\end{aligned}
$$
where we have used \eqref{u-f}. Consequently, 
$$
-\frac{\partial}{\partial k}[\rho-k\, f]= f
$$
and integrating this equation on $(0,1)$, we obtain
\begin{equation}
\label{RHO}
\rho(t, x, k)= k\, f(t, x, k) + \int_{k}^1 f\left(t, x, v\right) dv. 
\end{equation}
\begin{comment}
\begin{lemma}\label{Lem_bound} Let $f$ as defined above, then
$$
\begin{aligned}
f(t,x,k)\big(1-f(t,x,k)\big)\,\big|\,g^{\prime}(k)\,\big|\,&\leq \,f(t,x,k)\big(1-f(t,x,k)\big),
\\[5pt]
f(t,x,k)\big(1-f(t,x,k)\big)\,\big|\,g(k)\,\big|\,&\leq \,f(t,x,k)\big(1-f(t,x,k)\big),
\end{aligned}
$$
for a.e. $t \in(0, T)$, $x \in \Omega$, $k \in \mathbb{R}$.
\end{lemma}
\begin{proof}
Prova a grosso modo, primero recuerde que $g'(k)=1-2k$ ouseja nao eh limitado.

\medskip
Para mostrar the assertion faremos por caso.

CASO 1: $k<0$. Entao $f=1$ nada a mostrar.

CASO 2: $k>1$. Entao $f=0$ pronto.

CASO 3: $0\leq k \leq 1$. Entao $|g'(k)|\leq1$ logo segue o resultado.
\end{proof}
%%%%%%%
\begin{remark} From Remark \ref{rem_rho} and \eqref{Eq_f}, we arrive
$$
\rho(t, x, k)-u(t, x) f(t, x, k)=0
$$
for $k>1$ and $k<0$.
\end{remark}
\end{comment}
%%%%%%%%%%%%%%%
\begin{lemma}\label{lem_cota_R}
Let $\rho$ be given by \eqref{RHO}. Then, for any $T> 0$
$$
\big|R(t, x, k) \big| \,\leq  F(t,x,k)
$$
for almost everywhere  $t \in (0, T)$, $x \in \Omega$, and $k \in (0,1)$.
\end{lemma}
\begin{proof}
Follows from Lemma 3.1 in \cite{perth-dali-tams}. 
\end{proof}

%%%%%%%%%%%%%%%%
\medskip
Now, we are ready to show the strong convergence of the family $\{u_\varepsilon\}$,
which is the main issue to prove Theorem \ref{ETGS} (Main Theorem). 

\begin{theorem}
\label{FF}
We have $F= 0$ almost everywhere in $(0,\infty) \times \Omega \times \R$. 
\end{theorem}
\begin{proof}
First, let us recall that $R$ and $F$ are identically zero for $k \leq 0$ and $k \geq 1$.
Now let us define for $\delta> 0$ sufficiently small, the function $\zeta _{\delta }$ as follows
\begin{equation}
 \label{kzi}
\zeta _{\delta }(z):=\left\{
\begin{aligned}
&0,\hspace{20pt}\text{if }z< 0,
\\[5pt]
&\frac{z}{\delta }, \hspace{20pt} \text{if } 0 \leq z \leq \delta,
\\[5pt]
&1,\hspace{20pt} \text{if } z> \delta, \text{ } 
\end{aligned}
\right. 
\end{equation}%
and consider $\Psi_\delta(t,x,k)= \phi _{\delta}(t)\,\psi _{\delta}(x)\,\xi _{\delta}(k)$,
where 
$$
\begin{aligned}
\phi _{\delta}(t)&:=\zeta {_{\delta}}(t)-\zeta {_{\delta}}
(t-t_{0}+\delta ),\quad \text{ for } t_{0}\in (2\delta ,T)  ,
\\[5pt]
\psi _{\delta}(x)&:=\zeta _{\delta}(\mathrm{d}(x)),\hspace{2.6cm} \text{ for } x\in \Omega ,
\\[5pt]
\xi_{\delta}(k)&:=\zeta {_{\delta}}(k + \delta^{-1})-\zeta_{\delta}(k-\delta^{-1}),\quad \text{ for } k\in\R,
\end{aligned}
$$
where $\mathrm{d}(x)= \min_{y \in \Gamma }|x-y|$ is the distance function from $x \in \overline{\Omega }$
\ to $\ \Gamma$, and also $t_0$ is a Lebesgue point of the function 
$$
t \mapsto \int_{\Omega\times\R}F(t,x,k)\,dx\,dk.
$$
Then, choosing $\Psi_\delta(t,x,k)$ as a test function in the respective integral form of 
\eqref{eq_f(1-f)}, we get the inequality
$$
\begin{aligned}
-\iint_{\Omega_T} \int_{[0,1]} F\left\lbrace\,\partial_t \Psi_\delta(t,x,k)
+\,\mathbf{b}\cdot \nabla_{\!\!(k,x)}\Psi_\delta(t,x,k)\right\rbrace \,dx\,dk \ dt
\\[5pt]
+\iint_{\Omega_T}\int_{[0,1]} R(1-2f)\Psi_\delta(t,x,k) \,dx\,dk \ dt \leq 0.
\end{aligned}
$$
From the definition of the drift vector field $\mathbf{b}$, i.e. equation \eqref{DRIFTB}, it follows that
\begin{equation}
\label{FINAL}
\begin{aligned}
&\frac{1}{\delta}\int_{t_0-\delta}^{t_0} \iint_{\Omega \times [0,1]} F  \ \psi_{\delta}(x) \,dk\,dx\,dt
\leq \frac{1}{\delta} \int_{0}^{\delta} \iint_{\Omega \times [0,1]} F  \ \psi_{\delta}(x) \,dk\,dx\,dt
\\[5pt]
&+ \frac{1}{\delta}\int_{\{0\leq \mathrm{d}(x)\leq \delta\}}
\iint_{[0,T]\times [0,1]} F \,\phi_{\delta}(t) \,
   \,g'(k)\nabla\mathcal{K}_sc\cdot\nabla \mathrm{d}(x)\,dt\,dk\,dx
\\[5pt]
&+ \iint_{\Omega_T} \int_{[0,1]} F \ \phi_{\delta}(t) \ \psi_{\delta}(x) \ dk \,dx \ dt
\\[5pt]
&= I_1^\delta + I_2^\delta + \iint_{\Omega_T} \int_{[0,1]} F \ \phi_{\delta}(t) \ \psi_{\delta}(x) \ dk \,dx \ dt, 
\end{aligned}
\end{equation}
with the obvious notation and we have used Lemma \ref{lem_cota_R}. 

\medskip
Now, due to Proposition \ref{pro_intial} and applying the Dominated Convergence Theorem, we have 
\begin{equation}
\label{I1}
  \lim\sup_{\delta \to 0} I_1^\delta \leq \iint_{\Omega \times [0,1]} \big(f_0 - (f_0)^2 \big) \ dk \ dx= 0. 
\end{equation}
For the term $I_2^\delta$ one observes that, since $\ \Gamma$ is a
$C^2-$boundary, there exists a sufficiently small $\delta> 0$ such that, each point 
$x \in \Omega_{\delta}:=\left\{ {x} \in \Omega :\ \mathrm{d}(x)< \delta \right\}$ has a
unique projection $\mathbf{x}_{b}=\mathbf{x}_{b}(x)$ on the
boundary $\Gamma .$ 
For every $x \in \Omega_\delta$, we have 
$$
   \nabla \mathrm{d}(x)= -\nu(\mathbf{x}_{b}) + O(\delta )
$$ 
and the Jacobian of the change of variables 
$$
   \Omega_\delta \ni x \leftrightarrow (\mathbf{x}_{b},\tau) \in \Gamma \times (0,\delta)
   \quad \text{is equal to $\frac{D(x)}{D(\mathbf{x}_{b},\tau)}= 1 + O(\delta )$},
$$ 
where $\tau= \mathrm{d}(x)$. Therefore, we obtain 
\begin{equation}
\label{I2}
\begin{aligned}
   I_2^\delta &\leq \int_0^T \!\! \int_\Gamma \int_0^1 F^\delta \,\phi_{\delta}(t) \, |\nabla\mathcal{K}_sc \cdot \nu| \ dk \ d\mathbf{x}_{b} \ dt + O(\delta )
   \\[5pt]
   &= O(\delta ), 
\end{aligned} 
\end{equation}
where we have used that $\nabla \mathcal{K}_sc  \in L^{2}\left((0,T); \mathbf{H}^1(\Omega)\right)$, 
$\nabla \mathcal{K}_sc \cdot \nu= 0$ on $\Gamma$, 
$$
   F^\delta:= \frac{1}{\delta} \int_0^\delta F(\cdot,(\mathbf{x}_{b},\tau),\cdot) \ d\tau, 
$$
and $(\mathbf{x}_{b},\tau)$ forms an orthogonal coordinate system in a neighborhood of $\tau= 0$. 

\medskip
Finally, from \eqref{I1} and \eqref{I2} we get passing to the limit in \eqref{FINAL} as $\delta \to 0$,  
\begin{equation}
\iint_{\Omega \times [0,1]} F(t,x,k) \,dk\,dx \leq \int_0^t \iint_{\Omega \times[0,1]} F(t',x,k) \,dk\,dx \ dt' \label{eq_desig_F}
\end{equation}
for almost everywhere $t \in[0,T]$. Therefore, applying the Gronwall's Lemma, we obtain from \eqref{eq_desig_F}
$$
\iint_{\Omega_T} \int_{[0,1]}F(t,x,k) \ dk \ dx \ dt \leq 0,
$$
which implies the result.
\end{proof}

\bigskip
The above theorem implies that, the kinetic function $f$ takes only the values 0 and 1
almost everywhere in $(0,\infty) \times \Omega \times\R$,
and since $f$ is monotone decreasing on $k$, there exists a function $w= w(t,x)$, such that
$$
f(t,x,k)=\sgn^{+}(w(t,x)-k).
$$
Therefore for any $G\in C^{1}([0,1])$, with $G(0)=0$, it follows that 
\begin{equation*}
\ G\left( u_{\varepsilon }\right) =\int\limits_{0}^{1}G^{\prime
}(v)f_{\varepsilon }(\mathbf{\cdot },\mathbf{\cdot ,}v)\ dv\rightharpoonup
\int\limits_{0}^{1}G^{\prime }(v)f(\mathbf{\cdot },\mathbf{\cdot ,}v)\
dv= G(w)
\end{equation*}%
weakly star in $L^{\infty}(\Omega _{T})$, which implies $w= u$ 
almost everywhere, and
the strong convergence of the family $\left\{ u_{\varepsilon }\right\} $\vspace{1pt} to
$u$ \ in $L^{p}(\Omega _{T})$ for any $p<\infty .$
Then, we write \eqref{eq_entropy_perturbado} in weak sense using the entropy pair $\mathbf{F}(u, v)$, and jointly with 
the second equation in \eqref{eq_regulation}, also written in the weak sense, 
we derive passing to the limit as $\varepsilon \to 0$
that, the pair $(u,c)$ satisfies \eqref{CLE} and \eqref{EE}, which ends the proof of Theorem \ref{ETGS}. 
\section{Comments and Extensions}
\label{Comments}

One remarks that, the results established in this paper apply to some interesting 
correlated versions of the system \eqref{eq.main}. 

\medskip
1.  Let us consider for $s \in (0,1)$ the following system
\begin{equation}
\left \{
\begin{aligned}
    &\partial_t u + \dive \big(g(u) \,  \nabla \mathcal{K}_s c \big)= 0, \quad  \text{in $(0,\infty) \times \Omega$},
\\[5pt]
    &  (-\Delta_N+I_d)^{1-s} \ c= u,  \quad \text{in $\Omega$},
\\[5pt]
    & u|_{\{t=0\}}= u_{0},  \hspace{21pt} \text{in $\Omega$},
\\[5pt]
 & \nabla\mathcal{K}_s c\cdot \nu= 0, \hspace{22pt}  \text{on $\Gamma$},
\end{aligned}
\right.\label{eq_Delt+Id}
\end{equation}
where the operator $(-\Delta_N+I_d)^{1-s}$ is analogously defined by the spectral theory. 
Indeed, there exists a complete orthonormal 
basis $\{\vp_k\}^{\infty}_{k=0}$ of $L^2(\Omega)$, where
$\vp_k$ satisfies the following eigenvalue problem 
\begin{equation}
\left\lbrace
\begin{aligned}
(-\Delta+I_d) \vp_k&= \mu_k \, \vp_k,\quad\mbox{ in }\Omega,
\\[5pt]
\nabla \vp_k\cdot\nu&=0,\quad\quad\quad\mbox{ on }\Gamma.
\end{aligned}
\right.\nonumber
\end{equation}
Therefore, we have that $\vp_k$ is the eigenfunction corresponding to eigenvalue $\mu_k$,
which is given by $\mu_k=\lambda_k + 1$ for each $k\geq 0$, where the pair $(\vp_k,\lambda_k)$ is 
the solution of \eqref{Eq_Neumann}. 
Thus, applying the functional calculus we can define
\begin{equation}
\begin{aligned}
    (-\Delta_{N}+I_d)^{s} u:= \sum_{k= 0}^\infty (\lambda_k + 1)^{s} \, \langle u, \vp_k \rangle \ \vp_k. 
\end{aligned}\nonumber
\end{equation}

Now, we are alloyed to take $\mathcal{K}_s= \big((-\Delta_{N}+I_d)^{s}\big)^{-1}= (-\Delta_{N}+I_d)^{-s}$. 
Similar to the system \eqref{eq.main}, it is not difficult to show that the condition \eqref{eq_contabilidade} is satisfied, that is to say 
$$
   \int_\Omega c(t,x) \, dx= \int_\Omega u(t,x) \, dx.
$$

\medskip
2. Finally, let $0\leq \sigma\leq 1$ be fixed and consider for $s \in (0,1)$ the following system
\begin{equation}
\left \{
\begin{aligned}
    &\partial_t u + \dive \big(g(u) \,  \nabla \mathcal{K}_s c \big)= 0, \quad  \text{in $(0,\infty) \times \Omega$},
\\[5pt]
    &  (-\Delta_N + \sigma I_d)^{1-s} \ c + (1-\sigma) c= u,  \quad \text{in $\Omega$},
\\[5pt]
    & u|_{\{t=0\}}= u_{0},  \hspace{21pt} \text{in $\Omega$},
\\[5pt]
 & \nabla\mathcal{K}_s c\cdot \nu= 0, \hspace{22pt}  \text{on $\Gamma$}.
\end{aligned}
\right.\label{eq_Delt+sigmaId}
\end{equation}
Clearly, for $\sigma= 0$ we get the system \eqref{eq.main} and for $\sigma= 1$ we have \eqref{eq_Delt+Id}. 
For $0< \sigma< 1$  we may similarly define the operators 
$$
   (-\Delta_{N} + \sigma \, I_d)^{s} \quad \text{and} \quad  \mathcal{K}_s= (-\Delta_{N} + \sigma \, I_d)^{-s}. 
$$
This system does not satisfy exactly the condition \eqref{eq_contabilidade}. Indeed, we have
\begin{equation}
\label{NEWCOND}
  (\sigma^{1-s} + 1 - \sigma) \int_\Omega c(t,x) \, dx= \int_\Omega u(t,x) \, dx,
\end{equation}
and for each $\sigma \in (0,1)$, it follows that $(\sigma^{1-s} + 1 - \sigma) \in (0,1)$. 
Moreover, for any $\sigma \in [0,1]$ the system \eqref{eq_Delt+sigmaId} 
turns into \eqref{DAPER}  (at least formally) passing to the limit as $s \to 0^+$. 

%%%%%%%%%%%%%%%%%%%%%%%%
 \section*{Acknowledgements}
%%%%%%%%%%%%%%%%%%%%%%%%

Conflict of Interest: Author Wladimir Neves has received research grants from CNPq
through the grant  308064/2019-4, and also by FAPERJ 
(Cientista do Nosso Estado) through the grant E-26/201.139/2021. 

\medskip
Data sharing not applicable to this article as no datasets were 
generated or analysed during the current study.

%%%%%%%%%%%%%%

\bigskip

%\vspace{1cm}

\noindent \textsc{Wladimir Neves} \hfill \textsc{Gerardo Huaroto}\\
Instituto de Matem\'atica \hfill  Departamento de Matem\'atica \\
Universidade Federal do Rio de Janeiro\hfill Universidade Federal de Alagoas \\
Av. Athos da Silveira Ramos, 149  \hfill Av. Lourival Melo Mota, S/N
 \\
Rio de Janeiro, RJ, Brazil \hfill Maceio, Al, Brazil \\
CEP 21941-909 \hfill CEP 57072-970\\
\texttt{wladimir@im.ufrj.br} \hfill
\texttt{gerardo.cardenas@im.ufal.br}

\end{document}